\documentclass[10pt,article]{amsart}
\usepackage{curves}
\usepackage{amsmath}

\newtheorem {theorem}{Theorem}[section]
\newtheorem {proposition}[theorem]{Proposition}
\newtheorem {lemma}[theorem]{Lemma}
\newtheorem {assumption}[theorem]{Assumption}
\newtheorem {corollary}[theorem]{Corollary}
\newtheorem {definition}[theorem]{Definition}
\newtheorem {example}[theorem]{Example}
\newtheorem {remark}[theorem]{Remark}

\newcounter{conjecture}

\newcommand{\eqnsection}{
   \renewcommand{\theequation}{\thesection.\arabic{equation}}
   \makeatletter
   \csname @addtoreset\endcsname{equation}{section}
   \makeatother}

\newcommand{\be}{{\begin{equation}}}
\newcommand{\ee}{{\end{equation}}}
\def \bt{\begin{theorem}}
\def \et{\end{theorem}}
\def \bea{\begin{eqnarray}}
\def \eea{\end{eqnarray}}
\def \bas{\begin{eqnarray*}}
\def \eas{\end{eqnarray*}}

\def \bGa{\mathbf \Gamma}
\def \de{\delta}

\def \bDe{\mathbf \Delta}

\newcommand{\eps}{\varepsilon}

\newcommand{\won}{{\mbox{\bf 1}}}
\newcommand{\wzero}{{\mbox{\bf 0}}}

\def \om{\omega}

\def\mus{{\mu^0}}

\def \wh{\widehat}
\def \wt{\widetilde}

\newcommand{\ls}[1]
   {\dimen0=\fontdimen6\the\font \lineskip=#1\dimen0
\advance\lineskip.5\fontdimen5\the\font \advance\lineskip-\dimen0
\lineskiplimit=.9\lineskip \baselineskip=\lineskip
\advance\baselineskip\dimen0 \normallineskip\lineskip
\normallineskiplimit\lineskiplimit \normalbaselineskip\baselineskip
\ignorespaces }
\newcommand{\req}[1]{(\ref{#1})}

\def \bS{{\bf S}}
\def \bX{{\bf D}}
\def \bB{{\bf B}}
\def \bA{{\bf A}}
\def \bI{{\bf I}}
\def \bP{{\bf P}}
\def \bR{{\bf R}}
\def \bY{{\bf Y}}
\def \bZ{{\bf Z}}
\def \bx{{\bf x}}
\def \by{{\bf y}}

\def \BB{{\mathcal B}}
\def \CC{{\mathcal C}}
\def \BD{{D}}
\def \Cov{{\mathcal K}}
\def \CCb{{\mathcal D}(\bP_t)}
\def \CD{{\mathcal D}_f}
\def \CDD{\widehat{\mathcal D}_f}
\def \DD{{\mathcal D}}
\def \EE{{\mathcal E}}
\def \FF{{\mathcal F}}
\def \GGf{{\mathcal G}_f}
\def \GG{{\mathcal G}}
\def \HH{{\mathcal H}}

\def \LL{{\mathcal L}}

\def \sss{{\mathcal S}}

\def \({\left(}
\def \){\right)}

\def \Proof{{\bf Proof: }}

\def \bc{\begin{center} }
\def \ec{\end{center} }
\def\Bbb{\mathbb}

\begin{document}

\eqnsection
\newcommand{\reals}{{\Bbb{R}}}
\newcommand{\F}{{\mathcal F}}
\newcommand{\D}{{\mathcal D}}
\newcommand{\Fn}{{{\mathcal F}_n}}
\newcommand{\Gn}{{{\mathcal G}_n}}
\newcommand{\Hn}{{{\mathcal H}_n}}
\newcommand{\Fp}{{{\mathcal F}^p}}
\newcommand{\Gp}{{{\mathcal G}^p}}
\newcommand{\hm}{\HH^\varphi}
\newcommand{\beq}[1]{\begin{equation}\label{#1}}
\newcommand{\eeq}{\end{equation}}
\newcommand{\integers}{{\rm I\!N}}
\newcommand{\E}{{\Bbb E}}
\newcommand{\PPP}{{\Bbb P}}
\def\var{{\rm Var}}
\def\cov{{\rm Cov}}
\def\one{{\bf 1}}
\def\Rfg{R_{X^f,g}}
\def\leb{{\mathcal L}eb}
\def\Ho{{\mbox{\sf H\"older}}}  
\newcommand{\ffrac}[2]
  {\left( \frac{#1}{#2} \right)}
\newcommand{\calF}{{\mathcal F}}
\newcommand{\dfn}{\stackrel{\triangle}{=}}
\newcommand{\beqn}[1]{\begin{eqnarray}\label{#1}}
\newcommand{\eeqn}{\end{eqnarray}}
\newcommand{\oo}{\overline}
\newcommand{\uu}{\underline}
\newcommand{\Var}{{\rm \,Var\,}}
\def\squarebox#1{\hbox to #1{\hfill\vbox to #1{\vfill}}}
\renewcommand{\qed}{\hspace*{\fill}
            \vbox{\hrule\hbox{\vrule\squarebox{.667em}\vrule}\hrule}\smallskip}
\newcommand{\half}{\frac{1}{2}\:}
\newcommand{\beaa}{\begin{eqnarray*}}
\newcommand{\eeaa}{\end{eqnarray*}}
\newcommand{\calK}{{\mathcal K}}

\bibliographystyle{amsalpha}

\title[Fluctuation Dissipation Theorem]{
Markovian perturbation, Response and Fluctuation Dissipation Theorem}

\author[Amir Dembo\,\, Jean-Dominique  Deuschel] {Amir Dembo$^*$\,\,
Jean-Dominique Deuschel$^\dagger$}

\date{June 17, 2007; Revised: February 16, 2010.
\newline\indent
$^*$Research  partially supported by NSF grants 
\#DMS-0806211, \#DMS-0406042 and \#DMS-FRG-0244323.
\newline\indent
$^\dagger$Research partially supported by DFG grant \#663/2-3.
\newline
\newline
{\bf AMS (2000) Subject Classification:}
{Primary: 60J25, 82C05, Secondary: 82C31, 60J75, 60J60, 60K35}
}

\begin{abstract}
We consider the Fluctuation Dissipation Theorem (FDT) of statistical 
physics from a mathematical perspective. We formalize the 
concept of ``linear response function'' in the general framework 
of Markov processes. We show that for processes 
out of equilibrium it depends not only on the 
given Markov process $X(s)$ but also on the chosen perturbation of it. 
We characterize the set of all possible 
response functions for a given Markov process and show 
that at equilibrium they all satisfy the FDT. That is, if the
initial measure $\nu$ is invariant for the given Markov semi-group,
then for any pair of times $s<t$ and nice functions $f,g$, the
dissipation, that is, 
the derivative in $s$ of the covariance of $g(X(t))$ and $f(X(s))$ 
equals the infinitesimal response at time $t$ and direction $g$ 
to any Markovian perturbation that alters the 
invariant measure of $X(\cdot)$ in the direction of $f$ at time $s$.
The same applies in the so called FDT regime near equilibrium, i.e.
in the limit $s \to \infty$ with $t-s$ fixed, provided $X(s)$ converges
in law to an invariant measure for its dynamics.
We provide the response function of two generic Markovian perturbations 
which we then compare and contrast for 
pure jump processes on a discrete space, for finite dimensional 
diffusion processes,
and for stochastic spin systems.
\end{abstract}

\maketitle

\section{Introduction and outline}
One of the fundamental premises of 
statistical physics, the Fluctuation Dissipation Theorem (FDT),
follows from the assumption that the response of a system in 
thermodynamic equilibrium to a small external perturbation is 
the same as its relaxation after a spontaneous fluctuation. The FDT 
provides an explicit relationship between the 
equilibrium fluctuation properties of the thermodynamic 
system and its linear response (e.g. susceptibility), 
which involve out-of-equilibrium quantities.
As such it relates the dissipation of dynamics at thermal equilibrium 
of molecular scale (i.e. microscopic) models, with observable macroscopic 
response to external perturbations, allowing the use of microscopic 
models for predicting material properties 
(in the context of linear response theory). Among notable 
special cases of the FDT are the Einstein 
relation between particle 
diffusivity and its mobility \cite{einstein} (or the recent 
accounts in \cite{Hanney}, \cite{LR} and \cite{Lo}), 
and the Johnson-Nyquist formula \cite{nyquist} 
for the thermal noise in a resistor.

When deriving the FDT in equilibrium statistical physics 
one typically starts from a measure $\mus(\cdot)$ which is 
often a Gibbs measure, characterized by a Hamiltonian $H(\cdot)$, 
and a dynamics $X(s)$ for which this measure is 
invariant. For $\delta>0$ small, one perturbs the dynamics so 
it becomes invariant for the Gibbs measure corresponding to 
the Hamiltonian $H(\cdot)+\delta f(\cdot)$. 
The linear response measures the effect of applying such 
perturbation at time $s$ on the rate of change, as $\delta \downarrow 0$ 
in the value of a test function $g(\cdot)$ at time $t>s$ (often taking
$f=g$ to be the state variable of the dynamics in question).
This response function (of $s$ and $t$) is then compared 
to the rate of change in $s$ of the covariance
between $f(X(s))$ and $g(X(t))$ at equilibrium, 
in the non-perturbed dynamics, whereby the FDT states that the
ratio between these two functions is merely 
$\beta$, the inverse of the system's temperature. 


Whereas the FDT is well established and understood 
in physics, at least in or near thermal equilibrium, see \cite{kubo,kubo2}, 
our goal here is to provide its rigorous derivation from a mathematical
perspective, as a result about perturbations of Markovian semi-groups.  
This is easy to do in special (simple) cases, most notably, 
when dealing with a Markov process on a finite state space. Aiming 
here instead for a unified derivation, 
we formalize in Definition \ref{def:resp} the
concept of having a linear {\em response function} in the general framework
of a family of continuous time, homogeneous Markov processes
$X^f(\cdot)$ that are invariant for the measures 
$\mu^f(\cdot)=e^f \mus(\cdot)$ (see Assumption \ref{ass:1} 
for the precise setting). 
In our definition, the response function $\Rfg (s,t)$ depends on 
the initial position of the process at time $0$, 
thereby allowing us to study the effect of the initial measure
on the FDT relation. Though  
this function is uniquely defined per family $X^f(\cdot)$, 
it depends not only on the
given Markov process $X(\cdot)$ but also on the chosen perturbation 
$X^f(\cdot)$ of it (compare for example Theorems \ref{thm3a} and \ref{thm3b}).
In Theorem \ref{theo-fdt} we characterize the set of all possible
response functions for a given Markov process $X(\cdot)$ and show
that they all satisfy the FDT relation (\ref{eq:fdt}). 
It states that if the initial measure $\mus(\cdot)$ is invariant for the 
underlying Markov process $X(\cdot)$, then the dissipation, that is 
the derivative in $s$ of the covariance between 
$f(X(s))$ and $g(X(t))$ for $s<t$ equals the $\mus$-average of the
infinitesimal response $\Rfg (s,t)$ to any Markovian fluctuation 
$X^{\de f}(\cdot)$ that for $\de \downarrow 0$ alters the
invariant measure of $X(\cdot)$ in the direction of test function 
$f(\cdot)$ at time $s$, as registered at time $t$ via 
test function $g(\cdot)$. We show in Corollary \ref{cor-new1}
that this FDT relation holds in the limit $s \to \infty$
and $t-s$ fixed, whenever the initial measure is such that 
the law of $X(s)$ converges (in the appropriate sense) 
to an invariant measure $\mus(\cdot)$ and  
in Proposition \ref{prop5} we specify 
the set of all possible response functions in case of 
$\mu^f$-symmetric (i.e. reversible) Markovian perturbations.
Note that the FDT relation has to do with invariance of $\mus(\cdot)$
but does not require 
the Markov process $X(\cdot)$ to
be reversible with respect to $\mus(\cdot)$.  
To further demonstrate how widely this theory can be used, we
construct in Section \ref{sec:generic} 
two generic families of Markov perturbations 
satisfying Assumption \ref{ass:1}, which apply  
for {\em any} Markov process (subject only to mild restrictions
on the domain of certain generators). Namely, the {\em time change}
of Propositions \ref{prop1} and \ref{prop2a} and the {\em generalized 
Langevin dynamics} of Proposition \ref{prop2} (in the symmetric case)
and Proposition \ref{prop2c} (in the general, non-symmetric case).
The bulk of the mathematical work in this paper is in proving that 
these 
two generic families admit a response function per Definition 
\ref{def:resp}. This is done in Theorems \ref{thm3a} 
and \ref{thm3b} of Section \ref{sec:resp}
which also provide an explicit formula for 
the response function in each case. 

Moving to examples of specific Markov processes, 
in Proposition \ref{prop6}
we use the simple sufficient condition of 
Proposition \ref{prop0} for the existence 
of a response function, to show that essentially
all choices of the response function that are
possible per Theorem \ref{theo-fdt} are indeed attainable in the
context of pure jump processes on a discrete state space. We 
also demonstrate there how to create a host of perturbations
via cycle decomposition (for example, perturbation of 
a Metropolis dynamics, or of a Glauber dynamics).  
In Section \ref{sec:diff} we illustrate our  
generic Langevin perturbation for 
diffusion processes on connected, compact, smooth,
finite dimensional manifolds without boundary, showing 
in Proposition \ref{prop:diff} that it is generated in this case
by the addition of a smooth drift, which for symmetric diffusions
is of gradient form.
Finally, in Section \ref{sec:spin} we demonstrate the flexibility
of our framework by considering such a perturbation
for infinite dimensional diffusion processes associated with 
stochastic spin systems 
in the setting of Gibbs distributions. 

All our derivations and results apply even 
when the Markov process $X(\cdot)$ has more than one 
invariant measure (as is often the case in statistical 
physics when the system's temperature is sufficiently low),
and most of them apply also for non-reversible dynamics
(i.e. having non-symmetric semi-groups).

We note in passing that the appearance of $\beta$ in the FDT in
physics is merely due to the definition of Gibbs measure as
proportional to $e^{-\beta H}$, with $H$ the corresponding 
Hamiltonian, and not counting $\beta$ as part of the
perturbation $\delta f$. It makes more sense for the mathematical 
version of the FDT to not mention $\beta$ (so it in effect corresponds 
to doing statistical physics at $\beta=1$). Also, 
though from a mathematical point of view
the response function can not be defined solely in terms of the 
given Markov process $X(\cdot)$, this is never an issue in 
physics, whereby viewing the (Markov) dynamics as a classical
approximation of a quantum system, the perturbation of its
Hamiltonian in direction $f(\cdot)$ uniquely defines also 
the perturbed quantum system dynamics, hence its classical 
approximation $X^f(\cdot)$ (for example, see the physics based
derivation in \cite{ichiyanagi} of the FDT in quantum statistical 
mechanics). 


\section{General theory: FDT at or near equilibrium}

A continuous time, homogeneous,
strong Markov process $X(t)$ with values in a
complete, separable metric space $\sss$ is defined on a
fixed probability space $(\Omega, \FF, \PPP)$.
We assume that $X(t)$ has 
right continuous sample path and Markov semi-group 
$\bP_t h(x)=\E_x(h(X(t))) \in \CC_b(\sss)$ 
for any $h \in \CC_b(\sss)$. Let $\CCb$ denote the domain of $\bP_t$ 
with respect to the supremum norm. That is, 
the closed vector space 
$\CCb = \{ h \in \CC_b(\sss) : \|\bP_t h -h\|_\infty \to 0$ 
as $t \to 0 \}$ on which $\bP_t$ is strongly continuous and such 
that $\bP_t:\CCb \mapsto \CCb$ and 
denote by $\DD(\LL)$ the domain
of the generator $\LL:\DD(\LL) \mapsto \CCb$ of the semi-group $\bP_t$ 
(i.e. $\DD(\LL):= \{ h \in \CCb : t^{-1} (\bP_t h -h) \to \LL h$ 
in $(\CCb,\|\cdot\|_\infty)$ for
$t \downarrow 0 \}$), which is a dense subset of $(\CCb,\|\cdot\|_\infty)$. 
Recall that if 
$g\in \DD(\LL)$ then $t \mapsto \bP_t g: \reals_+ \mapsto \DD(\LL)$ is
differentiable and $\partial_t \bP_t g = \LL \bP_t g = \bP_t \LL g$.

For example, when $\sss$ is compact one often has $\CCb=\CC_b(\sss)$ 
whereas for $\sss=\reals^d$ one typically has 
$\CCb=\{h+c: c \in \reals, h \in \CC_0(\sss) \}$ 
for the space $\CC_0(\sss)$ of continuous functions that 
vanish at infinity.

The following hypothesis applies throughout. 
\begin{assumption}
We assume that $\CCb$ is an algebra under point-wise multiplication 
and a dense subset 
of $L_2(\mu)$ for any probability measure $\mu$. 
We consider test functions in a
linear vector space $\GG \subseteq \DD(\LL)$ 
such that $\GG$ is an algebra with
$\won \in \GG$ and that $\GG$ is dense in $(\CCb,\|\cdot\|_\infty)$. 
Setting $\wh{\GG} := \{ \bP_t g : g \in \GG, t \geq 0 \}$,
we further assume
that $\phi(g) \in \GG$, $\LL g \in \GG$ and $g \wh{h} \in \DD(\LL)$
for all $g \in \GG$, $\wh{h} \in \wh{\GG}$ and 
any $\phi \in \CC^\infty(\reals)$.
\end{assumption}

As usual we say that a finite measure $\mu$ on $\sss$ is
invariant for $\bP_t$ if $\langle \bP_t g \rangle_\mu :=
\int \bP_t g d\mu = \int g d\mu$ for all $g \in B_\sss$ and
$t \geq 0$. Since $\GG \subseteq \DD(\LL)$ is dense in
$L_2(\mu)$, this is equivalent to $\int \LL \wh{g} d\mu =0$
for all $g \in \wh{\GG}$. 

A key role is to be played by the symmetric bi-linear operator
$$
\bGa (f,g)=\LL(fg)-f\LL(g)-g\LL(f) \,,
$$ 
whose domain is
$\DD (\bGa) := \{ (f,g) : f \in \DD(\LL), g \in  \DD(\LL),
fg \in \DD(\LL) \} \subseteq \CCb \times \CCb$. One may  
instead define $\bGa$ as the carr\'e du champ operator
for the Dirichlet form associated with a $\mus$-symmetric 
semi-group $\bP_t$, thus possibly having a larger domain 
(c.f. \cite[Proposition I.4.1.3]{Bouleau}). 
However, in either case $\GG \times \GG$ is a 
dense subset of $\DD(\bGa)$, which is all we use in this paper. 

Let $\nu=\nu_0$ denote the initial measure of $X(0)$ and
$\nu_t := \nu 
\circ 
\bP_t$ the corresponding measure of $X(t)$.
Fixing $0 \leq s < t < \infty$ and $f \in L_2(\nu_s)$, $g \in L_2(\nu_t)$,
we denote the covariance of the
random variables $f(X(s))$ and $g(X(t))$
by $\Cov_{f,g}(s,t)$. By the Markov property we have,
\begin{equation}
\label{eq:cov0}
\Cov_{f,g}(s,t) = \langle \bP_s (f \bP_{t-s} g) \rangle_\nu -
\langle \bP_s f \rangle_\nu \langle \bP_t g \rangle_\nu \,.
\end{equation}

The next lemma provides useful formulas for the derivative of 
the covariance with respect to $s$.

\begin{lemma}\label{cov-lem}
For any $f,g \in \GG$ and all $0 \leq s < t < \infty$,
\begin{eqnarray}
\label{eq:cov1}
\partial_s \Cov_{f,g}(s,t)&=&
 \langle \bP_s  \LL (f \bP_{t-s} g) \rangle_\nu -
 \langle \bP_s (f \bP_{t-s} \LL g) \rangle_\nu -
\langle \bP_s \LL f \rangle_\nu \langle \bP_t g \rangle_\nu \\
&=& \int_\sss \bP_s \bGa ( f , \bP_{t-s}g ) d\nu + \Cov_{\LL f,g}(s,t) \,.
\label{eq:cov1b}
\end{eqnarray}
Suppose $\nu$ is invariant for $\bP_t$. Then,
\begin{equation}
\label{eq:cov2}
\partial_s \Cov_{f,g}(s,t) = -\langle f \LL \bP_{t-s} g \rangle_\nu \,.
\end{equation}
If in addition $\bP_t$ is $\nu$-symmetric
(i.e. $\langle f \bP_t g \rangle_\nu = \langle g \bP_t f \rangle_\nu$
for all $f,g \in B_\sss$),
then
\begin{equation}
\label{eq:cov3}
\partial_s \Cov_{f,g}(s,t)
= \frac{1}{2} \langle \bGa (f,\bP_{t-s} g) \rangle_\nu \,.
\end{equation}
\end{lemma}

\noindent
\Proof Fix $f,g$ and $s,t$ as in the statement of the lemma, and
$\de>0$. Let
$\BD$ denote the right-hand-side of (\ref{eq:cov1}), 
$h=\LL g$, $\bDe_\de=\bP_\de-\bI$ and $\wh{\bDe}_\de=\de^{-1} (\bP_\de-\bI)
-\LL$.
It is not hard to verify that
\beaa
&& \de^{-1}(\Cov_{f,g}(s+\de,t)-\Cov_{f,g}(s,t))  -\BD =
\langle \bP_s \wh{\bDe}_\de f \bP_{t-s} g \rangle_\nu
- \langle \bP_{s+\de} f \bP_{t-s-\de} \wh{\bDe}_\de g \rangle_\nu \\
&+& \langle \bP_{s+\de} f \bP_{t-s-\de} \bDe_\de h \rangle_\nu
- \langle \bP_s \bDe_\de f \bP_{t-s} h \rangle_\nu
+ \langle \bP_s \wh{\bDe}_\de f \rangle_\nu \langle \bP_t g \rangle_\nu \,.
\eeaa
Recall that $\{\bP_t\}$ is contractive for the supremum norm while 
$\| \bDe_\de h \|_\infty \to 0$ and
$\| \wh{\bDe}_\de h \|_\infty \to 0$ 
as $\de \downarrow 0$ for each fixed
$h \in \DD(\LL)$. Since $f$ is bounded with $f$, $g$, 
$f \bP_{t-s} g$ and $f \bP_{t-s} h$ in $\DD(\LL)$, it follows that 
$$
\lim_{\de \downarrow 0} \, 
\de^{-1}(\Cov_{f,g}(s+\de,t)-\Cov_{f,g}(s,t))  =\BD \,.
$$
Similar computation applies for the case of $\de<0$, thus establishing
\req{eq:cov1}.
Since $\bP_{t-s}\LL g = \LL \bP_{t-s} g$ for $g \in \GG$,
the equality \req{eq:cov1b} is a direct consequence of
\req{eq:cov0} and the definition of $\bGa (\cdot,\cdot)$.
To derive \req{eq:cov2} note that 
$$
\BD = \langle \LL (f \bP_{t-s} g) \rangle_{\nu_s} -
 \langle f \LL \bP_{t-s} g \rangle_{\nu_s} -
\langle \LL f \rangle_{\nu_s} \langle \bP_t g \rangle_\nu \,.
$$
If $\nu$ is invariant for $\bP_s$ then $\nu_s=\nu$ and further
$\langle \LL h \rangle_{\nu_s}=0$ for any 
$h \in \DD(\LL)$. Consequently, in this case 
$\BD=- \langle f \LL \bP_{t-s} g \rangle_\nu$ yielding \req{eq:cov2}.
If in addition $\bP_t$ is $\nu$-symmetric, then obviously
$\langle f \LL h \rangle_\nu =
\langle  h \LL f \rangle_\nu$ for all $f,h \in \DD(\LL)$, so
\req{eq:cov3} holds by the definition of $\bGa (\cdot,\cdot)$.
\qed

\begin{definition}\label{def-gconv}
The $\GG$-convergence as $s \to \infty$ of 
probability measures $\nu_s$ on $\sss$ 
to a probability measure $\mu$ on $\sss$, 
denoted $\nu_s \stackrel{\GG}{\to} \mu$, means that 
$\langle h \rangle_{\nu_s} \to \langle h \rangle_{\mu}$ as $s \to \infty$, 
for any fixed 
$h \in \{ f\wh{g}\,, \LL (f\wh{g})\, : f \in \GG,\, \wh{g} \in \wh{\GG} \}$.
\end{definition}
For example, the weak convergence of $\nu_s$ to $\mu$ implies
that $\nu_s \stackrel{\GG}{\to} \mu$ as well.

\begin{corollary}\label{cor0}
If $\nu_s =\nu \circ \bP_s \stackrel{\GG}{\to} \mu$ as $s \to \infty$,
then \req{eq:cov1} implies that for all $f,g \in \GG$, and any
fixed $\tau$,
$$
\lim_{s \to \infty} \partial_s \Cov_{f,g}(s,s+\tau)
=\langle \LL (f\bP_{\tau} g) \rangle_\mu -
 \langle f\bP_{\tau} \LL g \rangle_\mu -
\langle \LL f \rangle_\mu \langle\bP_\tau g \rangle_\mu \;.
$$
If in addition $\bP_t$ is $\mu$-symmetric
then by \req{eq:cov3},
$$
\lim_{s \to \infty} \partial_s \Cov_{f,g}(s,s+\tau)
= \frac{1}{2} \langle \bGa (f,\bP_{\tau} g) \rangle_\mu \,.
$$
\end{corollary}

Given an invariant probability measure $\mus$ for $\bP_t$ we
consider throughout the following setting (where 
the Banach space $\BB$ is typically $L_2(\mus)$ or 
$(\CC_b(\sss),\|\cdot\|_\infty)$). 
\begin{assumption}\label{ass:1}
The norm topology on a Banach space $(\BB,\|\cdot\|)$ of $\mus$-integrable 
functions is finer than the one induced by $L_1(\mus)$ and 
$\|h\| \le \|h\|_\infty$ for all $h \in \CCb \subseteq \BB$. 
For each $f \in \GG$ there exists a continuous time, homogeneous
Markov process $X^f(t)$ with a contractive semi-group $\bP_t^f$
on $(\BB,\|\cdot\|)$ such that 
$\|\bP_t^f h - h \| \to 0$ as $t \to 0$ for any $h \in \CCb$
and the finite, positive measure $\mu^f$ on $\sss$ such that
$d\mu^f/d\mus =e^f$ is invariant for $\bP_t^f$. 
Further, $X^0(t)=X(t)$, the subspace $\GG$ is contained in
the domain of the generator $\LL^f$ of $\bP_t^f$,
and $\LL^f g \in \GG$ for all $g \in \GG$.
\end{assumption}
\begin{remark} It is easy to verify that all results and proofs 
of this section remain valid in case $\mu^f = e^{\Phi(f)} \mus$ 
as long as $\Phi:\GG \to B_\sss$ is such that for any fixed
$f \in \GG$ the functions 
$\psi_\delta := \delta^{-1} (e^{\Phi(\delta f)}-\delta f -\won)$
converge to zero in $L_1(\mus)$ when $\de \downarrow 0$. 
\end{remark}

\medskip
The FDT is about the relation between derivatives of the covariance
at equilibrium and the {\em response} of the system to small
perturbation out of equilibrium. We turn now to the rigorous
definition of the latter (see also Proposition \ref{prop0} for 
an easy sufficient condition for existence of such response,
in case $\LL$ is bounded on $(\BB,\|\cdot\|)$).
\begin{definition}\label{def:resp}
Assume \ref{ass:1} and that for any
$f \in \GG$ there exists a linear operator $\bA_f: \DD(\bA_f) \mapsto \BB$
whose domain $\DD(\bA_f)$ contains $\wh{\GG}$,
and such that $s \mapsto \bA_f \bP_s g$ is strongly continuous
in $(\BB,\|\cdot\|)$ for each $g \in \GG$.
If moreover, for any $T \geq 0$, $\wh{g} \in \wh{\GG}$, 
all $t \in [0,T]$ and $\de>0$,
\begin{equation}\label{res:smpl}
\| \de^{-1} (\bP_t^{\de f} - \bP_t) \wh{g}
- \int_0^t \bP_s \bA_f \bP_{t-s} \wh{g} ds \| \leq \eta_\de t \,,
\end{equation}
and
$\eta_\de =\eta_\de(f,\wh{g},T) \to 0$ as $\de \downarrow 0$, then we call
\begin{equation}\label{res:form}
\Rfg (s,t) = \bP_s \bA_f \bP_{t-s} g \in \BB
\end{equation}
(which for any $t \geq 0$ is strongly integrable on $[0,t]$),
the {\em response function} 
(at time $t$ and direction $g$)
for the Markovian perturbations $X^f(\cdot)$ on $\GG$
(applied at time $s$).
\end{definition}

Note that $\Rfg (s,t)$ is uniquely defined per given
Markovian perturbations $X^f(\cdot)$. Indeed, suppose that
for some $f \in \GG$ the inequality
\req{res:smpl} holds for the same semi-groups $\bP_t^{\de f}$ and
both linear operators $\bA_f$ and $\wt{\bA}_f$.
Then, taking $\de \downarrow 0$ we see that
the linear operator $\bDe_A = \wt{\bA}_f - \bA_f$
is such that for all $t>0$ and $\wh{g} \in \wh{\GG}$,
\begin{equation}\label{eq:uniqr}
t^{-1} \int_0^t \bP_s \bDe_A \wh{g} ds +
t^{-1} \int_0^t \bP_{t-s} \bDe_A (\bP_s - \bI) \wh{g} ds = \wzero \,.
\end{equation}
With $\bP_s$ strongly continuous, upon taking $t \downarrow 0$ the
left-most term converges to $\bDe_A \wh{g}$, whereas
$t^{-1} \int_0^t \bP_{t-s} \bDe_A (\bP_s - \bI) \wh{g} ds \to \wzero$
by the contractivity of $\bP_{t-s}$ on $(\BB,\|\cdot\|)$ and the assumed
strong continuity of $\bDe_A \bP_s \wh{g}$.
Consequently, $\bA_f=\wt{\bA}_f$ on the set $\wh{\GG}$, and so using
either operator in (\ref{res:form}) leads to the same response function.

As we demonstrate next, $\Rfg (u,a+t+b)$ of (\ref{res:form}) is
merely the effect in ``direction'' $g$ and at time $a+t+b$,
of a small perturbation of the dynamics in ``direction'' $f$
during the time interval $u \in [a,a+t]$,
in agreement with the less formal definition of response function 
one often finds in the literature.
\begin{corollary}\label{cor1}
For any $f,g \in \GG$, each $T \geq a,b,t \geq 0$ and all $\de>0$,
a response function of the form (\ref{res:form}) must satisfy the
inequality
\begin{equation}\label{res:df}
\| \de^{-1} \bP_a (\bP_t^{\de f} - \bP_t) \bP_b g
- \int_a^{a+t} \Rfg (u,a+t+b) du \| \leq \eta_\de t \,,
\end{equation}
for $\eta_\de =\eta_\de(b,f,g,T) \to 0$ as $\de \downarrow 0$.
\end{corollary}

\noindent
\Proof Using the expression \req{res:form} to write \req{res:df}
more explicitly, one finds that for $a=0$ the latter is precisely
\req{res:smpl} for $v=b$, hence it obviously holds. Moreover,
in case $a>0$ we merely consider the norm of $\bP_a h^\delta$, where
$h^\delta$ is the element of $\BB$ the norm of which we consider
in \req{res:smpl}. As $\bP_a$ is contractive on $(\BB,\|\cdot\|)$, we
have that \req{res:df} holds also in this case.
\qed

We provide now an explicit sufficient condition for 
the existence of a response function of the form 
(\ref{res:form}) when $\LL^f$, $\LL$ and $\bA_f$ 
are bounded operators (as is the case in the
setting of Section \ref{sec:jump}). 
\begin{proposition}\label{prop0}
Assume \ref{ass:1} holds. If  
for each $f \in \GG$, the operators  
$\bA_f$, $\LL$, $\LL^{\de f}$, $\de >0$ are
bounded on $(\BB,\|\cdot\|)$ and the corresponding 
operator norms are such that
\begin{equation}\label{eq:am1}
\lim_{\de \downarrow 0} \| \de^{-1} (\LL^{\de f} - \LL) - \bA_f \| = 0 \,,
\end{equation}
then (\ref{res:smpl}) holds and $\bA_f \bP_s g$ is strongly continuous,
for each $g \in \GG$.
\end{proposition}

\noindent
\Proof Since $\|\bA_f \bP_{s+v} g - \bA_f \bP_v g \| \le \|\bA_f \|
 \, \|\bP_s \wh{g} - \wh{g} \|$
for $\wh{g} = \bP_v g$, the strong continuity of 
$s \mapsto 
\bA_f \bP_s g$ is a direct consequence
of the strong continuity of $\bP_s$ on its domain.
Fixing $f,g \in \GG$ and $v \geq 0$
it remains only to show that $t^{-1} \|\rho^\de_t \| \to 0$ 
as $\de \downarrow 0$, uniformly 
over $t \in (0,T]$, where 
$$
\rho^\de_t := \de^{-1} (\bP_t^{\de f} - \bP_t) \wh{g}
- \int_0^t \bP_{t-u} \bA_f \bP_{u} \wh{g} du \,. 
$$
To this end, let $r_t := - \int_0^t \bP_{t-u} \bA_f \bP_u \wh{g} du$
and note that 
$$
\partial_t \rho_t^\de =
\de^{-1} (\LL^{\de f}  \bP_t^{\de f} -\LL \bP_t) \wh{g} -
\bA_f \bP_t 
\wh{g}+\LL r_t \,,
$$
which imply, after some algebraic manipulations, that 
\begin{equation}\label{eq:deriv-rho}
\partial_t \rho^\de_t = \LL^{\de f} \rho^\de_t + \zeta^\de_t \,,
\end{equation}
for 
\begin{equation}\label{df:zeta}
\zeta^{\de}_t := (\LL - \LL^{\de f}) r_t +
[ \de^{-1} (\LL^{\de f} - \LL) - \bA_f ] \bP_t \wh{g} \,.
\end{equation}
It is not hard to show that the solution $\rho^\de_t$ 
of (\ref{eq:deriv-rho}) with initial condition $\rho^\de_0 = \wzero$, is
$$
\rho_t^\de = \int_0^t \bP^{\de f}_{t-s} \zeta^\de_s \,ds  \;.
$$ 
With both $\bP_t^{\de f}$ and $\bP_t$ contractive on $(\BB,\|\cdot\|)$,
we thus have that $\|r_t\| \leq t \|\bA_f\| \|g\|$ and for any $t \in [0,T]$,
$$
\| \rho_t^\de \| \leq \int_0^t  \| \zeta^\de_s \| ds \leq t [
 T \| \LL - \LL^{\de f}\| \| \bA_f \| + 
 \| \de^{-1} (\LL^{\de f} - \LL) - \bA_f \| ] \|g \| 
=: \eta_\de(f,g,T) t \,.
$$
Finally, note that from (\ref{eq:am1}) we have
that $\eta_\de \to 0$ as $\de \downarrow 0$. 
\qed

Our next theorem characterizes the type of response functions
one may find. It also proves the FDT, showing
that if $X^f(\cdot)$ has a response function $\Rfg (s,t)$, 
then the average of the response function according to an 
initial measure $\nu_0=\mus$ which is invariant for $X(\cdot)$, 
equals the time derivative of the covariance of $X(\cdot)$
under the same initial measure.
\begin{theorem}[Fluctuation Dissipation Theorem]\label{theo-fdt}
Let $f \in \GG$.
If $X^f(\cdot)$ has a response function, then $\bA_f \won = \wzero$ and
$\bA_{rf} = r \bA_f$ for all $r>0$. Further, then 
$\langle (\bA_f + f \LL) \wh{g} \rangle_{\mus} = 0$ for all
$\wh{g} \in \wh{\GG}$ and consequently, if the initial measure
$\nu_0=\mus$, then for any $s<t$,
\begin{equation}
\label{eq:fdt}
\partial_s \Cov_{f,g}(s,t) = \langle \Rfg (s,t) \rangle_{\nu_0} \,.
\end{equation}
\end{theorem}
\begin{remark}\label{rem:green-kubo} 
If $\bP_t$ is $\mus$-symmetric and 
$f,g \in L_2(\mus)$ with $(f,g) \in \DD(\bGa )$, then by spectral decomposition
we have the Green-Kubo formula 
$$
-\langle f \LL g \rangle_\mus = \frac{1}{2} \langle \bGa (f,g) \rangle_\mus
= \int_0^\infty [ \langle (\bP_s \LL f) (\LL g) \rangle_\mus ] ds    
$$
(c.f. \cite[Theorem 4.3.8]{JQQ}). Applying it for 
$f$ and $\bP_{t-s} g$, we get the alternative expression 
$$
\partial_s \Cov_{f,g}(s,t) =  
\int_0^\infty [ \langle (\bP_s \LL f) (\LL \bP_{t-s} g) \rangle_\mus ] ds
$$    
for the dissipation term.
In contrast with 
(\ref{eq:fdt}), this identity
does not involve a perturbation of the given Markovian dynamic.  
\end{remark}

\noindent
\Proof 
Fix $f \in \GG$.
If $g=\won$ then $\bP_u g = \bP_u^{\de f} g = \won$ for 
all $\de>0$ and $u \geq 0$, so in this case taking $\de \downarrow 0$ 
in (\ref{res:smpl}) we find that for any $t>0$,
$$
t^{-1} \int_0^t \bP_s \bA_f \won ds = \wzero \,.
$$
Thus, taking $t \downarrow 0$ we have that $\bA_f \won = \wzero$.

Next, fixing $r>0$, note that $\bP_t^{\de (rf)}= \bP_t^{(\de r)f}$ for all
$\de>0$, hence by (\ref{res:smpl}) we have for $\bDe_A=\bA_{rf}-r\bA_f$,
$$
\| \int_0^t \bP_s \bDe_A \bP_{t-s} \wh{g} ds \| \leq
[\eta_{\de} (rf,\wh{g},T)+r \eta_{\de r} (f,\wh{g},T)] t \,.
$$
Taking $\de \downarrow 0$ this implies that (\ref{eq:uniqr}) holds
for $\bDe_A=\bA_{rf}-r \bA_f$ and all $t>0$. Hence, by the argument 
we provided immediately following (\ref{eq:uniqr}),
we deduce that $\bDe_A \wh{g} = \wzero$ for all $\wh{g} \in \wh{\GG}$.
That is, without loss of generality we may assume that $\bA_{rf}=r \bA_f$, as
claimed.  

Let $\psi_\de = \de^{-1} (e^{\de f} - \de f - {\bf 1})$. 
Since $\mu^{\de f}$ is invariant for $\bP_t^{\de f}$ and $\mus$ is 
invariant for $\bP_t$, it follows that
\beaa
\langle \de^{-1} (\bP_t^{\de f} - \bP_t) \wh{g} \rangle_\mus
&=& \langle \psi_\de (\bI-\bP^{\de f}_t) \wh{g} \rangle_\mus
- \langle f  (\bP_t^{\de f} - \bP_t) \wh{g} \rangle_\mus
- \langle f (\bP_t - \bI) \wh{g} \rangle_\mus \\
&:=& F_1(\de) + F_2(\de) +F_3 \,.
\eeaa
With \req{res:smpl} implying that
$\| (\bP_t^{\de f} - \bP_t) \wh{g} \| \to 0$ as
$\de \downarrow 0$ 
(hence so does 
$\langle | (\bP_t^{\de f} - \bP_t) \wh{g} | \rangle_{\mus}$),
and $\| f\|_\infty <\infty$ we deduce that 
$F_2(\de) \to 0$. Further, $\langle | \psi_\de | \rangle_\mus \to 0$ 
when $\de \downarrow 0$ (since $f$ is bounded), 
while
$\| (\bI-\bP^{\de f}_t) \wh{g} \|_\infty 
 \leq 2 \| \wh{g} \|_\infty$ is uniformly
bounded in $\de$, resulting with $F_1(\de) \to 0$ as well. We thus deduce
by considering the limit $\de \downarrow 0$ in \req{res:smpl},
and applying Fubini's theorem, that for any $t>0$ and $\wh{g} \in \wh{\GG}$,
$$
\langle f t^{-1} (\bP_t - \bI) \wh{g} \rangle_\mus
+ t^{-1} \int_0^t \langle \bP_{t-s} \bA_f \bP_{s} \wh{g} \rangle_\mus ds = 0 \,.
$$
Further, with $\mus$ invariant for $\bP_t$ and having assumed strong
continuity of $
s \mapsto 
\bA_f \bP_s \wh{g}$, upon taking $t \downarrow 0$ we
find that
\begin{equation}\label{eq:just-done}
\langle f \LL \wh{g} \rangle_\mus
+ \langle \bA_f \wh{g} \rangle_\mus = 0 \,,
\end{equation}
as claimed.

Next, recall \req{eq:cov2} and \req{res:form} that for $\nu_0=\mus$
which is invariant for $\bP_s$ and for any finite $s < t$,
$$
\partial_s \Cov_{f,g}(s,t) - \langle \Rfg (s,t) \rangle_\mus
= -\langle f \LL \bP_{t-s} g \rangle_\mus - \langle  \bA_f \bP_{t-s} g \rangle_\mus
= 0 \,,
$$
where the right-most identity is precisely \req{eq:just-done}.
\qed

\medskip
Combining Corollary \ref{cor0} and Theorem \ref{theo-fdt} we 
deduce the existence of the FDT regime in out-of-equilibrium
dynamics whenever $X(s)$ converges in law (in the 
appropriate sense), with the limiting measure $\mus$ being 
invariant for $X(\cdot)$.
That is, the FDT relation (\ref{eq:fdt}) then asymptotically 
holds for $t-s=\tau$ fixed, in the limit $s \to \infty$. 
Specifically, similar to Definition \ref{def-gconv} we define
the notion of $\GGf$-convergence as follows.
\begin{definition}
The $\GGf$-convergence of probability measures $\nu_s$ on $\sss$ 
to a probability measure $\mu$ on $\sss$, 
denoted $\nu_s \stackrel{\GGf}{\to} \mu$, means that 
$\langle h \rangle_{\nu_s} \to \langle h \rangle_{\mu}$ as $s \to \infty$, 
for any fixed 
$h \in \{ f \wh{g},\, \LL (f\wh{g}), \, \bA_f \wh{g} \, 
: \wh{g} \in \wh{\GG} \}$ (implicitly 
assuming that $\bA_f \wh{g} \in L_1(\nu_s)$ for all $s$ large
enough).
\end{definition}
\begin{corollary}\label{cor-new1}
Let $\nu_0$ denote the initial measure of $X(0)$ 
for an $\sss$-valued, continuous time, 
homogeneous, strong Markov process $X(t)$ 
such that $X^f(\cdot)$ has a response function (in the
sense of Definition \ref{def:resp}). Suppose further that 
$\nu_s = \nu_0 \circ \bP_s \stackrel{\GGf}{\to} \mus$ as $s \to \infty$.
Then, for any $g \in \GG$ and fixed $\tau \geq 0$,
\begin{equation}\label{eq:qfdt}
\lim_{s \to \infty} \partial_s \Cov_{f,g}(s,s+\tau) =
- \langle f \LL \bP_\tau g \rangle_{\mus} 
= \lim_{s \to \infty} \langle \Rfg (s,s+\tau) \rangle_{\nu_0} 
\;.
\end{equation}
\end{corollary}

\noindent
\Proof The left side of the identity follows from the formula 
for the limit of $\partial_s \Cov_{f,g}(s,s+\tau)$ given
in Corollary \ref{cor0} and the fact that $\mus$ is invariant
for $\bP_t$ so $\langle \LL h \rangle_{\mus} = 0$ for all 
$h \in \DD(\LL)$. Recall (\ref{res:form}) that 
$\langle \Rfg (s,s+\tau) \rangle_{\nu_0} =
\langle \bA_f \bP_{\tau} g \rangle_{\nu_s}$ which 
converges to $\langle \bA_f \bP_{\tau} g \rangle_{\mus}$ 
as $s \to \infty$. From Theorem \ref{theo-fdt} we have that  
$\langle (\bA_f + f \LL) \bP_\tau g \rangle_{\mus} = 0$,
which thus yields the right side of the stated identity.  
\qed


\medskip
We conclude this section with a statement characterizing 
response functions for symmetric perturbations.
\begin{proposition}\label{prop5}
Assume \ref{ass:1} holds for the Hilbert space 
$\BB = L_2(\mus)$ and that for each $f \in \GG$
the semi-group $\bP^f_t$ of $X^f(t)$ is $\mu^f$-symmetric.
Then, any response function of the form (\ref{res:form})
is based on $\bA_f=\bB_f - f \LL$ for a linear operator $\bB_f$ 
which is $\mus$-symmetric on $\wh{\GG}$.
\end{proposition}
\begin{remark} The stated $\mus$-symmetry of $\bB_f = \bA_f+ f \LL$ 
is necessary in this setting of symmetric perturbations. However,
typically more is required from $\bB_f$ in order to assure the 
existence of a response function of the form (\ref{res:form}).
\end{remark}
\noindent
\Proof Note that for any $\delta>0$ and all 
$\wh{g}, \wh{h}\in\wh{\GG}$,  
\beaa
F_\delta (\wh{h},\wh{g}) &:=& 
\delta^{-1} \langle \wh{h} (\bP_t^{\delta f}-\bP_t)\wh{g}\rangle_\mus +
\langle \wh{h} f(\bP_t-\bI)\wh{g}\rangle_\mus\\
&=& \langle \wh{h} \psi_\delta(\bI-\bP_t^{\delta f})\wh{g}\rangle_\mus
-
\langle \wh{h} f(\bP_t^{\delta f}-\bP_t)\wh{g}\rangle_{\mus}\\
&+&\delta^{-1}\langle \wh{h} (\bP_t^{\delta f}-\bI)
\wh{g}\rangle_{\mu^{\delta f}}-
\delta^{-1}\langle \wh{h} (\bP_t-\bI)\wh{g}\rangle_{\mus}\\
&=:& F_1(\delta,\wh{h},\wh{g})+F_2(\de,\wh{h},\wh{g})+
F_3(\delta,\wh{h},\wh{g})+F_4(\delta,\wh{h},\wh{g}) \,,
\eeaa
where $\psi_\delta=\delta^{-1} (e^{\delta f} - \delta f - \won)$.
Recall that $\| (\bI-\bP_t^{\delta f})\wh{g} \|_\infty 
 \leq 2 \|\wh{g}\|_\infty$, $\wh{h}$ is bounded 
and $\langle |\psi_\delta| \rangle_\mus \to 0$, hence
$F_1(\delta,\wh{h},\wh{g})\to 0$ as $\delta\downarrow 0$.
Further, 
$\langle |(\bP_t^{\delta f}-\bP_t)\wh{g}| \rangle_{\mus} \to 0$
hence also $F_2(\delta,\wh{h},\wh{g})\to 0$ as $\delta\downarrow 0$.
By the $\mu^{\delta f}$-symmetry of 
$\bP^{\delta f}_t$ and the $\mus$-symmetry of $\bP_t$, it follows that    
$$F_3(\delta,\wh{h},\wh{g})=F_3(\delta,\wh{g},\wh{h}), \qquad
F_4(\delta,\wh{h},\wh{g})=F_4(\delta,\wh{g},\wh{h}), 
$$
for any $\delta>0$. Consequently, as $\delta \downarrow 0$,
$$
F_\delta (\wh{h},\wh{g}) - F_\delta(\wh{g},\wh{h})  \to 0
$$
which by (\ref{res:smpl}) and the $\mus$-symmetry
of $\bP_{t-u}$ amounts to 
\begin{equation}\label{eq:symmet}
E_t (\wh{h},\wh{g}) = E_t (\wh{g},\wh{h}) \,,
\end{equation}
where
\beaa 
E_t(\wh{h},\wh{g}) :=  
\int_0^t\langle (\bP_{t-s} \wh{h}) (\bA_f \bP_s\wh{g}) \rangle_{\mus}\,ds
+ \langle \wh{h} f(\bP_t-\bI)\wh{g}\rangle_{\mus} \,.
\eeaa
Since 
$s \mapsto 
\bA_f \bP_s \wh{g}$ is strongly continuous (as part of 
Definition \ref{def:resp} of the response function), it follows that
as $t \downarrow 0$,
$$
t^{-1} \int_0^t\langle \wh{h} \bA_f \bP_s\wh{g} \rangle_{\mus}\,ds
+ t^{-1} \langle \wh{h} f(\bP_t-\bI)\wh{g}\rangle_{\mus} \to 
\langle \wh{h} \bA_f \wh{g} \rangle_\mus +
\langle \wh{h} f \LL \wh{g}\rangle_{\mus}  =
\langle \wh{h} \bB_f \wh{g} \rangle_\mus \,.
$$
Further, for all $u>0$,
$$
(\bP_u-\bI)\wh{h}=\int_{0}^{u} \bP_v \LL\wh{h} \,dv  \,,
$$
hence $\| (\bP_u-\bI)\wh{h}\|_\infty \leq u\| \LL\wh{h}\|_\infty$,
while the strong continuity of $\bA_f \bP_s \wh{g}$ implies
that $\sup_{s \leq t} \| \bA_f \bP_s \wh{g} \| \to \|\bA_f \wh{g}\| <\infty$  
as $t \downarrow 0$. Taken together, these imply that   
$$ \lim_{t\downarrow 0}
t^{-1}\int_0^t\langle ((\bP_{t-s}-\bI)\wh{h}) 
(\bA_f \bP_s\wh{g})\rangle_{\mus}\,ds=0\,,
$$
and consequently we have also that 
$$
t^{-1} E_t(\wh{h},\wh{g}) \to \langle \wh{h} \bB_f \wh{g} \rangle_\mus \,.
$$ 
We thus conclude, based on (\ref{eq:symmet}), that 
$\langle \wh{h} \bB_f\wh{g}\rangle_{\mus}
=\langle \wh{g} \bB_f\wh{h}\rangle_{\mus}$
for any $\wh{g}, \wh{h} \in \wh{\GG}$, 
as claimed.
\qed

\section{Generic Markov perturbations}\label{sec:generic}
In this section we construct several generic Markov 
perturbations for which Assumption \ref{ass:1} holds.  
Our presentation is somewhat technical because we aim 
at addressing a rather general framework. The reader 
may thus benefit from considering first the 
concrete examples of Sections \ref{sec:jump} and \ref{sec:diff}.

\medskip
\noindent{\bf Time change}.
Our first construction corresponds to
changing the clock as follows. For each fixed $f \in \GG$ and $s \geq 0$ let
$$
t^f(s,\om):=\int_0^s e^{f(X(u))} du \,.
$$
Note that $t^f:\reals_+ \times \Omega \to \reals_+$ is a measurable
stochastic process the sample path of which are
everywhere differentiable with $\frac{d}{ds} t^f = e^{f(X(s))}$ bounded
and bounded away from zero. Its inverse,
$\tau^f(t) := \inf\{s \geq 0: t^f(s) \geq t\}$ is thus also a measurable
stochastic process, the sample path of which are everywhere differentiable
with $\frac{d}{dt} \tau^f = e^{-f(X(\tau^f(t)))}$ uniformly (in $t$ and
$\om$) bounded and bounded away from zero.
\begin{proposition}\label{prop1}
Assumption \ref{ass:1} holds for 
$(\BB,\|\cdot\|)=(\CC_b(\sss),\|\cdot\|_\infty)$
and the Markov process $X^f_0(t)=X(\tau^f(t))$. Further, the 
generator $\LL^f_0$ of the semi-group $(\bP^f_0)_t$ of $X^f_0(t)$ is 
such that $\DD(\LL) = \DD(\LL^f_0)$ and
$\LL^f_0 g =e^{-f} \LL g$ for any $g \in \DD(\LL)$.
\end{proposition}

\noindent
\Proof Obviously, $\{ \tau^f(t) \le s \} = \{t^f(s) \ge t\}$ is 
in $\FF_s = \sigma( X(u) : 0 \leq u \leq s)$ 
by the right continuity and boundedness of $u \mapsto \exp(f(X(u)))$. 
Hence, for each fixed $t$, the random variable $\tau^f(t)$ is
a stopping time with respect to the canonical filtration $\FF_s$
of $\{X(\cdot)\}$. Similarly, the stochastic process
$$
\tau^f (t) = \int_0^t e^{-f(X^f_0(v))} dv \,,
$$
is adapted to the canonical filtration
$\FF^f_t = \sigma( X^f_0 (u) : 0 \leq u \leq t)$ of $\{X^f_0 (\cdot)\}$.
The strict monotonicity and continuity of $t \mapsto \tau^f(t)$
imply that $\FF^f_t = \FF_{\tau^f(t)}$ and further it is not hard to check
that $\tau^f(t,\om)$ have the regeneration property
$$
\tau^f(t,\om)=\tau^f(s,\om)+\tau^f(t-s,\theta^{\tau^f(s,\om)} \om) \,,
$$
for any $t>s > 0$
(where $\theta^u \omega(\cdot) = \omega(u+\cdot)$ denotes the usual shift 
operator). 
Thus, for any $h \in B_\sss$,
$t>s >0$ and $x \in \sss$, by the strong Markov property of $X(\cdot)$ at
$\tau^f(s)$,
\bea
\nonumber
\E_x[h(X^f_0(t))|\FF^f_s] &=& \E_x[h(X(\tau^f(t)))|\FF_{\tau^f(s)}] \\
&=& \E_{X(\tau^f(s))} [h(\wt{X}(\tau^f(t-s))) ] =
\E_{X^f_0(s)} [h(\wt{X}^f_0(t-s)) ] \,,
\label{eq:mkv0}
\eea
where $\wt{X}(\cdot)$ and $\wt{X}^f_0(\cdot)$ denote independent copies
of $X(\cdot)$ and $X^f_0(\cdot)$, respectively. The Markov property of
$X^f_0(\cdot)$ then follows from \req{eq:mkv0} by standard arguments.
Since $X^f_0(\cdot)$ assumes its 
values in a complete, separable metric space,
the contractive semi-group $\bP_t^f h = \E_x(h(X^f_0(t)))$
is well defined on $(\CC_b(\sss),\|\cdot\|_\infty)$. 
Further, by the change of variable $v=\tau^f(u)$ 
its resolvent is given by
$$
\bR_\lambda^f h = \E_x \Big[ \int_0^\infty e^{-\lambda u} h(X^f_0(u)) du  \Big]
= \E_x \Big[ \int_0^\infty e^{-\lambda t^f(v)} e^{f(X(v))} h(X(v)) dv \Big] \,.
$$
Note that $\lambda t^f(s)= \lambda c s - \int_0^s \xi(X(u)) du$ 
for positive, finite $c=\exp(\|f\|_\infty)$ and 
the continuous function 
$\xi(x)=\lambda (c-e^{f(x)}) \ge 0$ with $\|\xi\|_\infty < \lambda c$. 
The linear operator $\bR_{\lambda c} \xi:\CCb \mapsto \DD(\LL)$ 
such that 
$$
(\bR_{\lambda c} \xi) g  = \E_x 
\Big[ \int_0^\infty e^{-\lambda c s} \xi(X(s)) g(X(s)) ds \Big]
$$
is thus strictly contractive, and since the series 
$$
\bR_\lambda^f h = \sum_{k \ge 0} (\bR_{\lambda c} \xi)^k 
\bR_{\lambda c} (e^f h)
$$
converges uniformly, it follows that 
$\bR_\lambda^f : \CCb \mapsto \CCb$
and consequently $\bP_t^f : \CC_b(\sss) \mapsto \CC_b(\sss)$.

Fixing $g \in \DD(\LL)$
recall that $g(X(\cdot))$ and $(\LL g) (X(\cdot))$ are
bounded right-continuous functions, with
$M(t):=g(X(t))-g(X(0))-\int_0^t (\LL g) (X(v)) dv$ a right-continuous
martingale with respect to the filtration $\FF_t$. With $M(0)=0$ and
$\tau^f(t)$ bounded above uniformly in $\om$, by Doob's optional
sampling theorem $\E_x M(\tau^f(t))=0$ for any $x \in \sss$.
By the change of variable $v=\tau^f(u)$ and Fubini, this amounts to
\begin{equation}\label{eq:Lf0}
\bP_t^f g(x) = g(x) + \int_0^t \E_x [ (e^{-f} \LL g) (X^f(u))] du
= g(x) + \int_0^t (\bP^f_u e^{-f} \LL g) (x) du \,.
\end{equation}
Let $\wt{g}=e^{-f} \LL g \in \CCb$. Since $\bP_u^f$ is contractive 
on $\CCb$,
it follows from \req{eq:Lf0} that
$\|\bP_t^f g -g \|_\infty \leq t \|\wt{g}\|_\infty$.
In particular, $\| \bP_t^f g - g \|_\infty \to 0$ as $t \downarrow 0$,
for all $g \in \GG$. With $\GG$ dense in the closed vector space
$(\CCb,\|\cdot\|_\infty)$, we have
that $\bP_t^f$ is strongly continuous there. The uniformly bounded
and strongly continuous $\bP^f_u \wt{g} : [0,t] \mapsto \CC_b(\sss)$
is also strongly integrable, so \req{eq:Lf0} implies that
for any $g \in \DD(\LL)$,
$$
\lim_{t \downarrow 0} \| t^{-1} (\bP_t^f g - g) - \wt{g} \|_\infty
= \lim_{t \downarrow 0}
\| t^{-1} \int_0^t \bP^f_u \wt{g} du - \wt{g} \|_\infty = 0\,.
$$
We have thus seen that $\DD(\LL) \subseteq \DD(\LL^f_0)$ with $\LL_0^f g =
e^{-f} \LL g$ for all $g \in \DD(\LL)$.

We next prove that if $g \in \DD(\LL^f_0)$, then necessarily
$g \in \DD(\LL)$. To this end, observe that
the sample path of $X^f(t)$ inherits the right continuity of those of
$X(t)$, and fixing $g \in \DD(\LL^f_0)$ (which is thus continuous 
and bounded), we have the right-continuous martingale
$M^f(s)=g(X^f(s))-g(X^f(0))-\int_0^s (\LL^f_0 g) (X^f(v)) dv$
with respect to the filtration $\FF^f_s$. With $M^f(0)=0$ and
the stopping time $t^f(s)$ for the latter filtration
bounded above uniformly (in $\om$), by the optional
sampling theorem $\E_x M^f(t^f(s))=0$ for all $x \in \sss$.
Since $X^f(t^f(u))=X(u)$, using Fubini and
the change of variable $v=t^f(u)$ we get in analogy to \req{eq:Lf0} that
$$
\bP_s g(x) = g(x)
+ \int_0^s (\bP_u e^f \LL^f_0 g) (x) du  \;,
$$
which by the uniform boundedness and strong continuity of
$\bP_u e^f \LL^f_0 g$ results with $g \in \DD(\LL)$
(and $\LL g = e^f \LL^f_0 g$), completing the proof that
$\LL_0^f = e^{-f} \LL$.

In particular, $\GG \subseteq \DD(\LL^f_0)$ and
since $\GG$ is an algebra containing both $e^{-f}$ and
$\LL g$, it follows that $\LL_0^f g \in \GG$ for any $f,g \in \GG$.
It remains just to verify that $\mu^f$ is invariant for $\bP_t^f$,
that is, $\int \LL^f_0 \bP_t^f g d\mu^f =0$ for all $t>0$ and $g \in \GG$.
We have already shown that
$\int \LL^f_0 \bP_t^f g d\mu^f = \int \LL (\bP_t^f g) d\mus$.
The latter is zero since
$\bP_t^f g \in \DD(\LL^f_0)=\DD(\LL)$ and $\mus$ is invariant for
$\bP_t$, so the proof of the proposition is complete.
\qed

\medskip
\noindent{\bf Perturbations for symmetric processes}. 
Relying on the powerful technology of Dirichlet forms (c.f. \cite{Bouleau}),
the second generic case we consider is that of a generalized
Langevin dynamics for a $\mus$-symmetric process.  

\begin{proposition}\label{prop2}
Suppose $\bP_t$ is $\mus$-symmetric. Then, for any $f \in \GG$ 
there exists a Markov process $X^f_1$ 
such that the generator of its $\mu^f$-symmetric, strongly 
continuous semi-group $(\oo{\bP}^f_1)_t$ 
on the Hilbert space $\BB=L_2(\mu^f)$, satisfies
\begin{equation}\label{eq:L1df}
\oo{\LL}^f_1 g =\LL g + \frac{1}{2} e^{-f} \bGa (e^f,g)\,, \qquad
\forall g \in \GG
\end{equation}
and for which Assumption \ref{ass:1} holds.
\end{proposition}

The Leibniz rule $\bGa (fh,g)=h\bGa (f,g)+f \bGa (h,g)$ applies 
whenever $\LL$ is the generator of a Markov 
process of continuous trajectories $t \mapsto X(t)$ 
(c.f. \cite{bakry}), resulting with  
$\oo{\LL}^f_1g=\LL g +\frac{1}{2}\bGa (f,g)$ as soon as
$\bGa(\sum_{k>n} \frac{f^k}{k!},g) \to 0$ when $n \to \infty$.
We call $X^f_1$ a generalized Langevin dynamics since
for $\sss$ a finite dimensional, compact, connected 
smooth manifold without boundary, 
the perturbed process $X^f_1(\cdot)$ is then
obtained by adding to the original (diffusion) 
process a drift of a gradient form (c.f. Section \ref{sec:diff}).

\begin{remark}
Proposition \ref{prop2} provides us with a
semi-group $(\oo{\bP}^f_1)_t$ of $X^f_1$ that is defined
only for $\mus$ almost every $x \in \sss$. In most interesting
specific situations
one easily shows that $(\oo{\bP}^f_1)_t$ is the unique extension 
to $L_2(\mu^f)$ of
a semi-group $(\bP^f_1)_t$ that is
strongly continuous on $(\CCb,\|\cdot\|_\infty)$
such that $(\bP_1^f)_t h (x) = \E_x( h(X^f_1(t)))$ for {\it all} $x \in \sss$
and $h \in \CC_b(\sss)$. Our proof of the 
proposition also shows that if $f(\cdot)$ is constant on $\sss$ then
$\oo{\LL}_1^f=\oo{\LL}$ is merely the closure of $\LL$ for 
the Hilbert space $L_2(\mus)$.
\end{remark}

\noindent
\Proof Fixing $f \in \GG$, we construct the continuous time,
homogeneous Markov process $X^f_1$ on the Hilbert space
$\HH=L_2(\mu^f)$.
To this end, consider the bi-linear form 
\begin{equation}\label{eq:efdf}
\EE_f(h,g) := \langle h (-\LL^f_1) g \rangle_{\mu^f}  
\,, \qquad \forall g,h \in \CD\,,
\end{equation}
where $\CD := \{g: e^f g \in \DD(\LL)\} \cap \DD(\LL)$ 
is a linear subspace of $\HH$
and $\LL_1^f g :=\LL g + \frac{1}{2} e^{-f} \bGa (e^f,g)$ 
is a linear operator from $\CD$ to $\HH$.
Our assumptions imply that the algebra $\GG$ is a
subspace of $\CD$, hence the latter is dense in $\HH$.
Recall that $\bP_t$ is $\mus$-symmetric, so the
same applies for its generator $\LL$ (i.e.
$\langle h \LL g \rangle_\mus = 
\langle g \LL h \rangle_\mus$ for all $g,h \in \DD(\LL)$).
It is easy to check 
that thus $\LL_1^f$ is $\mu^f$-symmetric 
and consequently, $\EE_f(\cdot,\cdot)$
is a symmetric form on $\CD \times \CD$.
Further, by the definition of $\LL$ and 
the $\mus$-symmetry of $\bP_t$ we find that for any $g \in \CD$,
\begin{equation}\label{eq:gtef}
\EE_f(g,g)= 
\frac{1}{2} \langle g^2 \LL e^f \rangle_{\mus}
- \frac{1}{2} \langle g e^f \LL g \rangle_{\mus}
- \frac{1}{2} \langle g \LL (e^f g) \rangle_{\mus} 
= \lim_{t \to 0} \EE_{f,t}(g)  \,,
\end{equation}
where for any $t>0$,
\begin{equation}\label{eq:etdf}
\EE_{f,t}(g)  := \frac{1}{2t} 
\langle g^2 - 2 g \bP_t g + \bP_t g^2 \rangle_{\mu^f} 
= \frac{1}{2t} \langle \E_x [ (g(X(t))-g(x))^2 ] \rangle_{\mu^f} \geq 0 \,,
\end{equation}
and consequently $\EE_f(g,g) \geq 0$ as well. 
The non-negative quadratic form $\EE_f$ on the dense subspace 
$\DD[\EE_f]=\CD$, is then closable (c.f. \cite[Example I.1.3.4]{Bouleau}). 
The closure $\oo{\EE}_f$ of $\EE_f$ determines a unique 
strongly continuous semi-group $(\oo{\bP}^f_1)_t$ 
of self-adjoint contractions on $\HH$ (c.f. \cite[Proposition
I.1.2.3]{Bouleau}). The generator $\oo{\LL}^f_1$ 
of $(\oo{\bP}^f_1)_t$ is (up to a sign inversion) the  
Friedrichs extension of $-\LL^f_1$, that is, a non-positive 
self-adjoint operator on $\HH$ satisfying \req{eq:L1df} 
for all $g \in \CD$ (c.f. \cite[Example I.1.3.4]{Bouleau}).

For $\epsilon>0$ let $\varphi_\epsilon: \reals \to [-\epsilon,1+\epsilon]$
be infinitely differentiable functions, such that 
$\varphi_\epsilon(t)=t$ for $t \in [0,1]$ and
$0 \leq \varphi_\epsilon(t_2)-\varphi_\epsilon(t_1) \leq t_2-t_1$
if $t_1 \leq t_2$ (see \cite[Problem I.1.2.1]{Fuk} for a
construction of such functions). Obviously, for any $\epsilon,t >0$,
$g \in \CCb$ and $\mu^f$ almost every $x \in \sss$, 
$$
\E_x[ (\varphi_\epsilon(g(X(t)) - \varphi_\epsilon(g(x)))^2]
\leq \E_x [ (g(X(t))-g(x))^2 ] \,.
$$
By (\ref{eq:etdf})
this implies that $\EE_{f,t}(\varphi_\epsilon(g)) \leq \EE_{f,t}(g)$. 
Further, recall that when 
$g \in \GG$ also $\varphi_\epsilon(g) \in \GG \subseteq \CD$,
in which case it follows from (\ref{eq:gtef}) that 
$\EE_f(\varphi_\epsilon(g),\varphi_\epsilon(g)) \leq \EE_f(g,g)$. 
As this holds for all $\epsilon>0$ and on the dense subset $\GG$
of $\DD[\oo{\EE}_f]$, we conclude that 
$\oo{\EE}_f$ is a symmetric Dirichlet form
in $\HH$ (c.f. \cite[Proposition I.4.10]{MR}). 
Consequently, the strongly continuous semi-group
$(\oo{\bP}^f_1)_t$ of self-adjoint contractions on $\HH$ is sub-Markovian 
(c.f. \cite[Proposition I.3.2.1]{Bouleau}), 
and in particular $\mus$-almost everywhere $(\oo{\bP}^f_1)_t \psi \geq 0$ 
whenever $\mus$-almost everywhere $\psi \geq 0$ (c.f.
\cite[Definitions I.2.1.1 and I.2.4.1]{Bouleau}). Recall that
$\won \in \GG$ and $\LL \won  = \wzero$ (since $\bP_t$ is Markovian),
implying by \req{eq:L1df} that $\oo{\LL}^f_1 \won = \wzero$ as well.
Consequently, $(\oo{\bP}^f_1)_t = e^{t \oo{\LL}^f_1}$ (c.f. \cite[Proposition 
I.1.2.1]{Bouleau}), is Markovian as claimed 
(that is, also $(\oo{\bP}^f_1)_t \won = \won$ for all $t>0$).
\qed

In carrying out the construction of Langevin dynamics in 
a non-symmetric setting we use the following analog of
Proposition \ref{prop1} which applies in the Hilbert 
space $L_2(\mu^f)$ for any $\mus$-symmetric Markov semi-group 
$\oo{\bP}_t$ on $L_2(\mus)$ (possibly no longer defined point-wise,
and even when the corresponding Markov process 
is neither strong Markov, nor has right continuous sample path). 
\begin{proposition}\label{prop2a}
If 
$\oo{\bP}_t$ is $\mus$-symmetric Markov semi-group 
on the Hilbert space $L_2(\mus)$ with a generator
$\oo{\LL}$, then Assumption \ref{ass:1} holds for a 
$\mu^f$-symmetric, strongly continuous Markov 
semi-group $(\oo{\bP}^f_0)_t$ on $\BB=L_2(\mu^f)$, 
whose generator $\oo{\LL}^f_0$ has the same domain as $\oo{\LL}$ and
is such that 
\begin{equation}\label{eq:L0df}
\oo{\LL}^f_0 g =e^{-f} \oo{\LL} g \,, \qquad \forall g \in \DD(\oo{\LL}) \,.
\end{equation}
\end{proposition}

\noindent
\Proof Since $\oo{\bP}_t$ is a $\mus$-symmetric 
sub-Markovian semi-group on $L_2(\mus)$, with domain that
is dense in $L_2(\mus)$ (on account of $\GG \subseteq \DD(\oo{\LL})$
being dense in $L_2(\mus)$), its generator $\oo{\LL}$ is 
a Dirichlet operator, namely, 
a (negative), self-adjoint 
operator on $L_2(\mus)$
of a dense domain on which $\langle (h-1)_+ \oo{\LL} h \rangle_{\mus} \leq 0$
(c.f. \cite[Proposition I.3.2.1]{Bouleau}). 
Fixing $f \in \GG$, the linear map $\oo{\LL}_0^f$  from the dense subset
$\DD(\oo{\LL})$ of $L_2(\mu^f)$ to $L_2(\mu^f)$ given by
(\ref{eq:L0df}) is a symmetric operator on $L_2(\mu^f)$, as
$$
\langle h_1 \oo{\LL}^f_0 h_2 \rangle_{\mu^f} =
\langle h_1 \oo{\LL} h_2 \rangle_{\mus} =
\langle h_2 \oo{\LL} h_1 \rangle_{\mus} =
\langle h_2 \oo{\LL}^f_0 h_1 \rangle_{\mu^f} 
$$
for all $h_1,h_2 \in \DD(\oo{\LL})$. Further,
if $h,\phi \in L_2(\mu^f)$ are such that
$\langle h \oo{\LL}^f_0 g \rangle_{\mu^f} = \langle \phi g \rangle_{\mu^f}$
for all $g \in \DD(\oo{\LL}^f_0)=\DD(\oo{\LL})$ then of course
$\langle h \oo{\LL} g \rangle_{\mus} = \langle e^f \phi g \rangle_{\mus}$
for $h,\phi \in L_2(\mus)$. With $\oo{\LL}$ self-adjoint on
$L_2(\mus)$ this implies that $h \in \DD(\oo{\LL})$ and
$\phi=e^{-f} \oo{\LL} h = \oo{\LL}^f_0 h$, namely, the symmetric operator
$\oo{\LL}^f_0$ is self-adjoint on $L_2(\mu^f)$.
Clearly, 
$$
\langle (h-1)_+ \oo{\LL}^f_0 h \rangle_{\mu^f} 
= \langle (h-1)_+ \oo{\LL} h \rangle_{\mus} \leq 0
$$
for all $h \in \DD(\oo{\LL})=\DD(\oo{\LL}^f_0)$. 
We thus deduce that $\oo{\LL}^f_0$
is a Dirichlet operator on $L_2(\mu^f)$, hence the generator of 
a $\mu^f$-symmetric sub-Markovian semi-group on $L_2(\mu^f)$, 
denoted hereafter $(\oo{\bP}_0^f)_t$  (for example, 
apply again 
\cite[Proposition I.3.2.1]{Bouleau}, now in the converse direction). 
Finally, with $\oo{\bP}_t$ a Markov semi-group, we have that
${\bf 1} \in \DD(\oo{\LL})$ and $\mus$-a.s.  
$\oo{\LL} {\bf 1} = {\bf 0}$. Of course, the same applies
$\mu^f$-a.s. for  $\oo{\LL}^f_0$ of (\ref{eq:L0df}). Consequently,
$(\oo{\bP}^f_0)_t$ is actually a $\mu^f$-symmetric 
Markov semi-group on $L_2(\mu^f)$ (e.g. \cite[Exercise I.3.1]{Bouleau}),
as claimed.
\qed


\medskip
\noindent{\bf Langevin dynamics in a non-symmetric setting}. 
Building on Propositions \ref{prop1} and \ref{prop2} 
we shall construct a 
generalized Langevin dynamics when $\mus$ is invariant for the 
strongly continuous Markov semi-group $\bP_t$ on $\CCb$ 
(of a strong Markov process 
of right continuous sample path $X(t)$ valued on a complete, separable metric 
space $\sss$), now in case $\bP_t$ is not $\mus$-symmetric. 
This construction is more complex than what we
have seen in Proposition \ref{prop2} and it involves
certain (mild) restrictions on the domains of various generators. 
We start by letting $\oo{\bP}_t$  
denote the extension of $\bP_t$ to a $\mus$-invariant 
strongly continuous, Markov semi-group on the Hilbert space $L_2(\mus)$
and $\oo{\bP}_t^\star$ the adjoint of $\oo{\bP}_t$ in $L_2(\mus)$.
The adjoint semi-group $(\oo{\bP}_t^\star)_{t \geq 0}$ is then 
strongly continuous in $L_2(\mus)$ and its
generator $\oo{\LL}^\star$ is the adjoint of the generator $\oo{\LL}$ of 
$\oo{\bP}_t$ 
(see \cite[Theorem 4.3]{goldstein}). 
Further, since
$\oo{\bP}_t$ is $\mus$-invariant and Markovian on $L_2(\mus)$, 
the same applies for $\oo{\bP}_t^\star$.

Let $\oo{\bGa}(h,g)=\oo{\LL} (gh) - h \oo{\LL}(g) - g \oo{\LL} (h)$
denote the $L_2(\mus)$ extension of $\bGa$. 
We seek a perturbed Markovian semi-group whose generator has the form
$$
\oo{\LL}^f_1 = 
\wt{\LL} + \frac{1}{2} e^{-f} \oo{\bGa}(e^f,\cdot) + e^{-f} \LL_a \,,
$$
where $\wt{\LL}=\frac{1}{2}(\oo{\LL}+\oo{\LL}^\star)$ and
$\LL_a =\frac{1}{2}(\oo{\LL}-\oo{\LL}^\star)$ 
correspond to the symmetric and anti-symmetric parts of $\oo{\LL}$,
respectively. Indeed, it is easy to see that $\oo{\LL}^f_1$ 
is $\mu^f$-invariant and $\oo{\LL}_1^{\de f} \to \LL$ 
for $\de \downarrow 0$. However, while $\wt{\LL}$ is a generator of a
Markovian semi-group (as shown in Corollary \ref{cor2b}), this is not
the case for $\LL_a$. Hence, we alternatively construct 
$\oo{\LL}^f_1$ in Proposition \ref{prop2c} as the sum 
of the Markovian generators $\oo{\LL}^f_0$ of 
Proposition \ref{prop2a} and $\wt{\LL}^f$ 
which we construct next. 
\begin{lemma}\label{lem2b}
Suppose $\bP_t$ is a $\mus$-invariant strongly continuous 
Markov semi-group. Assume further that 
$\GG \subseteq \DD(\oo{\LL}^\star)$ and $\oo{\LL}^\star g \in \GG$
for all $g \in \GG$. Then, for each non-negative $f \in \GG$ there
exists a Markov process $\wt{X}^f$ such that the generator 
of its $\mu^f$-symmetric, strongly continuous semi-group
$\wt{\bP}^f_t$ on the Hilbert space $\BB=L_2(\mu^f)$, satisfies
\begin{equation}\label{eq:wLdf}
\wt{\LL}^f g = (1-e^{-f}) \Big(\frac{1}{2} \oo{\LL} g
             +\frac{1}{2} \oo{\LL}^\star g \Big) + \frac{1}{2} e^{-f}
\oo{\bGa}(e^f,g) \,,
\qquad \forall g \in \GG
\end{equation}
and for which Assumption \ref{ass:1} holds.
\end{lemma}

\noindent
\Proof Fixing a non-negative $f \in \GG$, consider the linear subspace 
\begin{equation}\label{eq:whddef}
\CDD := 
\{ h : e^f h \in \DD(\oo{\LL}) \} \cap \DD(\oo{\LL}) \cap \DD(\oo{\LL}^\star)
\end{equation}
of $\HH=L_2(\mu^f)$ and the linear operator $\wh{\LL}^f$ from  
$\CDD$ to $\HH$, defined via (\ref{eq:wLdf}). 
It is not hard to verify that
$\wh{\LL}^f$ is $\mu^f$-symmetric operator on $\CDD$, with the
associated symmetric bi-linear form 
\bea\label{eq:ewtdf}
\wh{\EE}_f (h,g) &:=& - \langle h \wh{\LL}^f g \rangle_{\mu^f} \\
&=& \frac{1}{2} \Big[ \langle h \oo{\LL} (g - g e^f) \rangle_\mus
+ \langle g \oo{\LL} (h - h e^f) \rangle_\mus 
+ \langle g h \oo{\LL} e^f \rangle_\mus \Big] \,.
\nonumber
\eea
Since $\oo{\LL}$ and its adjoint $\oo{\LL}^\star$ are the generators of the 
$\mus$-invariant, strongly continuous semi-groups
$\oo{\bP}_t$ and $\oo{\bP}^\star_t$ on $L_2(\mus)$,
similarly to (\ref{eq:gtef}), we further have that
for any $g \in \CDD$,
\begin{equation}\label{eq:gtwef}
\wh{\EE}_f(g,g)= 
\frac{1}{2} \langle g^2 \oo{\LL} e^f \rangle_{\mus}
+ \langle g \oo{\LL} g \rangle_{\mus}
- \langle g \oo{\LL} (e^f g) \rangle_{\mus} 
= \lim_{t \to 0} \wh{\EE}_{f,t}(g)  \,,
\end{equation}
where for any $t>0$,
\begin{equation}\label{eq:wetdf}
\wh{\EE}_{f,t}(g)  := \frac{1}{2t} 
\langle (e^f-1) [g^2 - 2 g \oo{\bP}^\star_t g + \oo{\bP}^\star_t g^2 ] 
\rangle_{\mus} 
\geq 0 \,,
\end{equation}
in view of the non-negativity of $f$. Hence,
$\wh{\EE}_f(g,g) \geq 0$ for all $g \in \CDD$, and since
by our assumptions $\GG \subseteq \CDD$, we conclude
that the non-negative quadratic form $\wh{\EE}_f$ of dense domain
$\CDD$ is closable, denoting by $\wt{\EE}_f$ its closure, 
and by $\wt{\bP}_t$ and $\wt{\LL}^f$ the corresponding strongly 
continuous semi-group and generator. Since 
$\oo{\bP}^\star_t$ is a Markov semi-group,
by the same
argument as in the proof of Proposition \ref{prop2} we deduce
that replacing $g$ by $\varphi_\epsilon(g)$ 
reduces the value of $g^2 - 2 g \oo{\bP}^\star_t g + \oo{\bP}^\star_t g^2 
\geq 0$,
implying that $\wh{\EE}_{f,t}(\varphi_\epsilon(g)) \leq \wh{\EE}_{f,t}(g)$
for all $\epsilon,t>0$ and $g \in \GG$, a dense subset of
$\DD[\wt{\EE}_f]$. Arguing again as in Proposition \ref{prop2},
we conclude that $\wt{\EE}_f$ is a symmetric Dirichlet form in $\HH$ 
with $\wt{\bP}_t=e^{t \wt{\LL}^f}$ a strongly continuous 
Markovian semi-group (recall that 
$\won \in \GG$ and $\oo{\LL} \won = \oo{\LL}^\star \won = \wzero$). 
\qed

\begin{remark} Though Assumption \ref{ass:1} holds for 
the $\mu^f$-symmetric Markov semi-group $\wt{\bP}^f_t$,
the Fluctuation Dissipation Theorem \ref{theo-fdt} fails 
in this case, for 
$\wt{\LL}^{\de f}$ of \eqref{eq:wLdf} converges to {\em zero} 
as $\de \downarrow 0$ (and not to $\LL$). We thus need the additional
ingredients of Proposition \ref{prop2c}
in order to complete the construction
of the Langevin dynamics in non-symmetric setting.
\end{remark}

Considering Lemma \ref{lem2b} for $f=\won$ and then 
Proposition \ref{prop2a} for $f=\log(1-e^{-1}) \won$ 
we get the existence
of the following symmetric Markov process.
\begin{corollary}\label{cor2b}
For any $\mus$-invariant strongly continuous Markov semi-group
$\bP_t$, if 
$\GG \subseteq \DD(\oo{\LL}^\star)$ 
and $\oo{\LL}^\star g \in \GG$
for all $g \in \GG$, then  
Assumption \ref{ass:1} holds for a Markov process 
$\wt{X}$ on the Hilbert space $L_2(\mus)$. The generator
of its $\mus$-symmetric, strongly continuous semi-group
$\wt{\bP}_t$ 
is such that for any $g \in \DD(\oo{\LL}) \cap \DD(\oo{\LL}^\star)$,
\begin{equation}\label{eq:wLd}
\wt{\LL} g = \frac{1}{2} \oo{\LL} g
             +\frac{1}{2} \oo{\LL}^\star g  \,.
\end{equation}
\end{corollary}

We complete our construction by adding to the 
generator of Lemma \ref{lem2b} an appropriate 
non-symmetric perturbation (taken from Proposition \ref{prop1}). 
\begin{proposition}\label{prop2c}
Suppose as in Lemma \ref{lem2b}
that $f \in \GG$ is non-negative and  
$\oo{\LL}^\star (\GG) \subseteq \GG$.
Further assume that $\DD(\oo{\LL})$ is 
contained in the domain of the
generator $\wt{\LL}^f$ of $\wt{X}^f$.
Then, Assumption \ref{ass:1} 
holds for a $\mu^f$-invariant Markov process $X^f_1$
on $\BB=L_2(\mu^f)$. 
The generator $\oo{\LL}^f_1$ of its strongly continuous 
semi-group $(\oo{\bP}^f_1)_t$ has 
the same domain as $\oo{\LL}$ and is such that 
\begin{equation}\label{eq:Ladf}
\oo{\LL}_1^{f} = \wt{\LL}^{f}  + e^{-f} \oo{\LL} \,.
\end{equation}
\end{proposition}

\begin{remark} Proposition \ref{prop2c}
is the non-symmetric generalization of Proposition \ref{prop2}
since if $\bP_t$ is $\mus$-symmetric then $\oo{\LL}^\star=\oo{\LL}$ and 
the operator of (\ref{eq:Ladf}) and the corresponding 
Markov process coincide with $\oo{\LL}^f_1$ (and $X^f_1$)
of Proposition \ref{prop2}. Also note that 
since each $f \in \GG$ is bounded below,   
you have the non-negativity condition of the proposition 
by adding to any given $f \in \GG$ 
a sufficiently large constant, without changing the 
corresponding normalized invariant measure $\mu^f/\mu^f(\sss)$.
\end{remark}

\noindent
\Proof From Proposition \ref{prop1}
we get the Markov process $X_0^f(t)$ on $(\CC_b(\sss),\|\cdot\|_\infty)$
whose semi-group is generated by $\LL_0^f = e^{-f} \LL$ and is 
invariant for $\mu^f$.
Let $\oo{\bP}_0^f$ denote 
the extension of this semi-group to $L_2(\mu^f)$
and $\oo{\LL}_0^{f}$ its generator. Since $f$ is bounded it
is easy to see that $\oo{\LL}_0^f = e^{-f} \oo{\LL}$ has
the same domain as $\oo{\LL}$. 
%
Subject to the assumptions of Lemma \ref{lem2b} 
we get the $\mu^f$-symmetric Markov semi-group $\wt{\bP}^f_t$
on the same Hilbert space, whose generator $\wt{\LL}^f$ satisfies
(\ref{eq:wLdf}). In particular, by Hille-Yosida
theorem we know that $\wt{\LL}^f$ is a 
closed operator, 
so our assumption that its domain
contains the domain of $\oo{\LL}$ implies that 
the operator $\oo{\LL}_1^f$ (on $\DD(\oo{\LL})$),
given by (\ref{eq:Ladf}) is a generator 
of a strongly continuous semi-group 
$\oo{\bP}^f_1$ (c.f. \cite[pages 631-639]{DS}). 
Further, with $e^f \in \GG$ (an
algebra) such that both $\LL (\GG) \subseteq \GG$
and $\oo{\LL}^\star (\GG) \subseteq \GG$, we have the same
properties for $\oo{\LL}^f_0$ and $\wt{\LL}^f$, hence 
for $\oo{\LL}_1^f$ as well. The latter operator 
is the sum of  
the generators of two $\mu^{f}$-invariant Markov semi-groups, 
so by Trotter's product formula 
(c.f. \cite[Theorem 8.12]{goldstein}), 
we conclude that the corresponding generated semi-group is
both Markovian and $\mu^{f}$-invariant.
\qed

\section{The corresponding response functions}\label{sec:resp}

The response function is not determined just by the given
Markov process $X(\cdot)$, its invariant measure $\mus$,
and the perturbation function $f \in \GG$. Indeed, 
Theorem \ref{thm3a} provides the response function
of the Markovian perturbations of
Propositions \ref{prop1} and \ref{prop2a}, 
while Theorem \ref{thm3b} provides the  
response function for the Markovian perturbation of
Proposition \ref{prop2c} under the uniform control 
of \eqref{eq:sector}
on the anti-symmetric part of the generator $\oo{\LL}$.
See also Corollary \ref{cor3c} for the simpler
response function in the symmetric case of Proposition \ref{prop2}. 
Typically the response functions in these two theorems 
are not the same. 
\begin{theorem}\label{thm3a}
Taking $(\BB,\|\cdot\|)=(\CC_b(\sss),\|\cdot\|_\infty)$,
the Markov perturbation $X^f_0(t)$ has a response function
of the form (\ref{res:form}), with $\bA_f = (\bA_0)_f = -f \LL$.
\end{theorem}

\noindent
\Proof Recall that $X^f_0(t)$ satisfies Assumption \ref{ass:1} 
for $(\BB,\|\cdot\|)=(\CC_b(\sss),\|\cdot\|_\infty)$. Fixing $f \in \GG$,
observe that 
$\DD((\bA_0)_f)=\DD(-f \LL)=\DD(\LL)$ which contains $\wh{\GG}$
by our assumptions about $\GG$. Further, by strong continuity of 
$\bP_s$ on $(\CCb,\|\cdot\|_\infty)$ 
it follows that for each given $g \in \GG$,
$$
(\bA_0)_f \bP_s g = -f \LL \bP_s g = -f \bP_s (\LL g) \,,
$$ 
is also strongly continuous. Fixing $\wh{g}=\bP_v g$ for
$v \geq 0$ and $g \in \GG$, 
and taking $\bP_t^f = (\bP_0^f)_t$,   
as in the proof of Proposition \ref{prop0}
it suffices to show that $\|\zeta^\de_t \|_\infty \to 0$
as $\de \downarrow 0$, uniformly 
over $t \in (0,T]$, with $\zeta^\de_t$ given by 
(\ref{df:zeta}). To this end, since  
$\LL \bP_u = \bP_u \LL$, $\LL_0^f= e^{-f} \LL$ and $(\bA_0)_f=-f \LL$, it 
is not hard to verify that 
$\zeta^\de_t =  \phi_\de \LL r_t - \psi_{-\de} \bP_t \wh{h}$,
where $\wh{h}=\bP_v h$ for $h=\LL g \in \GG$, the continuous functions 
$\phi_\de := \won-e^{-\de f}$ and
$\psi_\de := \de^{-1} (e^{\de f} - \de f -\won)$ are such that 
$\| \phi_\de \|_\infty + \| \psi_{-\de} \|_\infty \to 0$ as $\de \downarrow 0$,
and $r_t = \int_0^t \bP_{t-u} f \bP_u \wh{h} du$
is in $\DD(\LL)$ by our assumptions about the space $\GG$ of test functions.

With $\bP_t$ contractive on $(\CCb,\|\cdot\|_\infty)$,
we thus have that for any $t \in [0,T]$,
$$
\| \zeta_t^\de \|_\infty \leq 
\| \phi_\de \|_\infty \| \LL r_t \|_\infty + 
\| \psi_{-\de} \|_\infty \|h \|_\infty \leq \eta_\de (f,g,T) \,,
$$
with $\eta_\de \to 0$ as $\de \downarrow 0$, provided 
$\sup \{ \| \LL r_s \|_\infty  : s \leq T \} < \infty$.

To show the latter,  recall that    
\beaa 
\LL r_s + f \bP_s \wh{h} &=& \partial_s r_s 
        = \partial_s \big( \int_0^s \bP_{u} f \bP_{s-u} \wh{h} \, du \big) \\
       &=&  \bP_s f \wh{h} + 
\int_0^s \bP_{u} f  \LL \bP_{s-u} \wh{h} \, du 
= \bP_s f \wh{h} + \int_0^s \bP_u f \bP_{s-u} \LL \wh{h}
\,du \,,    
\eeaa
hence by contractivity of the semi-group $\bP_u$,
\beaa
\| \LL r_s \|_\infty  &\leq& \| f \bP_s \wh{h} \|_\infty 
+ \| \bP_s f \wh{h} \|_\infty 
+ \int_0^s \| \bP_u f \bP_{s-u} \LL\wh{h} \|_\infty \, du 
\\
 &\leq& 2 \| f \|_\infty \| h \|_\infty + 
T \| f \|_\infty \| \LL h \|_\infty  \,,
\eeaa
which completes the proof.
\qed

\begin{theorem}\label{thm3b}
Suppose in addition to the assumptions of Proposition \ref{prop2c}
holding for $\delta f$ and all $\delta \in (0,1]$, that 
$(\oo{\bP}^{\de f}_1)_t \wh{g} \in \DD(\oo{\LL}^\star)$
for all $\wh{g} \in \wh{\GG}$, $t \geq 0$, and that 
$\DD(\oo{\LL}) \subseteq \DD(\wt{\LL})$ for $\wt{\LL}$ of 
Corollary \ref{cor2b}, with 
\begin{equation}\label{eq:sector}
\|\wt{\LL} r -\oo{\LL} r\|^2 
\leq - K \langle r \oo{\LL} r \rangle_\mus \,, 
\end{equation}
for some $K<\infty$ and all $r \in \DD(\oo{\LL})$.
The semi-groups $(\oo{\bP}^f_1)_t$ for the Markovian perturbations
$X^f_1 (t)$ on $\BB=L_2(\mus)$ then have a response function
of the form (\ref{res:form}), now with
\begin{equation}\label{eq:aalpha}
\bA_f = (\bA_1)_f =
\frac{1}{2} 
\oo{\bGa} (f,\cdot)+ f \wt{\LL}  - f \oo{\LL} \,,
\end{equation}
a linear operator on $\DD((\bA_1)_f) := 
\{ g : f g \in \DD(\oo{\LL}) \} \cap \DD(\oo{\LL})$.
\end{theorem}

\begin{remark} Since $(\bA_1)_f$ is applied in (\ref{res:form})
only on functions in the set 
$\wh{\GG}=\{\bP_t g : g \in \GG, t \geq 0 \}$, it is not hard to 
verify that (\ref{eq:aalpha}) is valid even with 
$\oo{\LL}$ and $\oo{\bGa}$ replaced by $\LL$ and $\bGa$, 
and that only the term corresponding to $f \wt{\LL}$ in (\ref{res:form})
might not be in $\CC_b(\sss)$ (hence requiring the application 
of $\oo{\bP}_s$ instead of $\bP_s$).
\end{remark}

\begin{remark}\label{rem:inconv}
Of the conditions of Theorem \ref{thm3b}, it is the least convenient
to check that $r=(\oo{\bP}^{\de f}_1)_t \wh{g}$ is in the domain 
of $\oo{\LL}^\star$. However, we use it only once, to deduce that then
$\langle r \oo{\LL} r \rangle_\mus = \langle r \wt{\LL} r \rangle_\mus$.
Thus, one can eliminate the former condition whenever there
is a direct way to verify the latter.
\end{remark}
Heuristically, the condition (\ref{eq:sector}) tells us that  
the symmetric part of the dynamics is dominant, in the same 
spirit as our assumption that $\DD(\oo{\LL})$ is contained
in both $\DD(\wt{\LL})$ and $\DD(\wt{\LL}^f)$ (and many
non-symmetric examples have different 
domains for $\oo{\LL}$, $\oo{\LL}^\star$ and $\wt{\LL}$).
Nevertheless, in many cases we do not have to worry about domains
of the various operators and in certain settings arrive at
the same conclusion even when (\ref{eq:sector}) does not hold (see
Section \ref{sec:jump} for one such example). 

\begin{remark}
Since $\mu^{\de \won} = \mus$ for any $\de>0$, it is 
natural to choose $\bP_t^{\de \won} = \bP_t$, as is the case for $X_1^f(\cdot)$,
yielding then that $\bA_{\won} = \wzero$.
However, this is not always done. 
For example, $\bP_t^{\de \won} \neq \bP_t$ 
for $X_0^f(\cdot)$, where indeed $(\bA_0)_\won = -\LL \neq \wzero$.
Also note that by Definition \ref{def:resp}
the response function $\Rfg (s,t)$ is always 
linear in $g$ and homogeneous with respect to multiplication of $f$ by
a positive scalar. In 
Theorems \ref{thm3a} and \ref{thm3b}, 
the response function is further linear in 
$f$, but this does not always apply
(see
Example \ref{exmp1}
for the Metropolis perturbation in which $\bA_{-f} \neq - \bA_f$). 
\end{remark}

\smallskip
In case of $\mus$-symmetric processes we get the following corollary 
upon considering $\wt{\LL}=\oo{\LL}^\star=\oo{\LL}$. 
The direct proof of this corollary is of course simpler and shorter 
than that of Theorem \ref{thm3b}.
\begin{corollary}\label{cor3c}
Suppose $\bP_t$ is $\mus$-symmetric and 
that for all $\delta \in (0,1]$ the generator of the
$\mu^{\de f}$-symmetric semi-group of Proposition \ref{prop2}
is such that $\DD(\oo{\LL}) \subseteq \DD(\oo{\LL}^{\de f}_1)$.
These semi-groups then have a response function
of the form (\ref{res:form}) with
$(\bA_1)_f = \frac{1}{2} \bGa(f,\cdot)$
(whose domain 
$\{ g \in \DD(\LL) : f g \in \DD(\LL) \}$ 
contains $\wh{\GG}$).
\end{corollary}

\noindent{\bf Proof of Theorem \ref{thm3b}:} 
Recall Proposition \ref{prop2c} that fixing hereafter 
non-negative $f \in \GG$, Assumption \ref{ass:1} holds
for the Hilbert space $\HH=L_2(\mus)$ and further,
$\wh{\GG}=\{\bP_t g : g \in \GG, t \geq 0 \}$ is
in the domain $\DD(\bA)$ of the linear operator $\bA=(\bA_1)_f$ 
of (\ref{eq:aalpha}).
Due to the linearity of $\bA$, for
the strong continuity of $\bA \bP_s g$ in $\HH$ per given $g \in \GG$,
it suffices to show the convergence to zero in $\HH$ of   
$\oo{\LL} (\bP_s \wh{g} - \wh{g})$, 
$\bGa (f, \bP_s \wh{g} - \wh{g})$ and
$\wt{\LL} (\bP_s \wh{g} - \wh{g})$, 
when $s \downarrow 0$, per given $\wh{g} \in \wh{\GG}$.
To this end, note first that 
$\oo{\LL} (\bP_s \wh{g} - \wh{g}) = (\bP_s-\bI) \LL \wh{g}$ 
(since $\wh{g}, \bP_s \wh{g} \in \DD(\LL)$), with the
latter converging to zero in $\HH$ by the strong continuity of 
$\bP_s$. Further, applying (\ref{eq:sector}) for 
$\bP_s \wh{g} -\wh{g} \in \DD(\oo{\LL})$, by 
the preceding argument, the convergence of 
$\wt{\LL} (\bP_s \wh{g} - \wh{g})$ in $\HH$ is a consequence of 
the convergence to zero of 
$$
- \langle (\bP_s \wh{g} - 
\wh{g}) \oo{\LL} (\bP_s \wh{g} - \wh{g}) \rangle_{\mus}
=
- \langle (\bP_s \wh{g} - \wh{g}) (\bP_s -\bI) (\LL \wh{g}) \rangle_{\mus}
\leq 
 2 \|\bP_s \wh{g} - \wh{g} \|  \, \|\LL \wh{g} \| \,,
$$
by the strong continuity of $\bP_s$. 

To deal with the last remaining term, namely 
$\bGa (f, \bP_s \wh{g} - \wh{g})$, note that
if $(h_1,h_2) \in \DD(\bGa)$ then 
$\bGa_t(h_1,h_2) \to \bGa(h_1,h_2)$ in supremum norm 
and hence also in $\HH$, when $t \downarrow 0$,  where 
\begin{equation}\label{eq:gtdf}
\bGa_t(g,h):=t^{-1} 
[g h - g \bP_t h - h \bP_t g + \bP_t g h ] \,,
\end{equation}
is a bi-linear symmetric, non-negative definite operator.

\noindent
Since
$\bGa_t(h_1,h_2)^2 \leq \bGa_t(h_1,h_1) \bGa_t(h_2,h_2)$, 
we have in particular that
\bea
\| \bGa (f, \bP_s \wh{g} - \wh{g}) \|^2
&=& \lim_{t \to 0} \| \bGa_t(f,\bP_s \wh{g} - \wh{g}) \|^2 
\nonumber \\
&\leq& \lim_{t \to 0} \| \bGa_t(f,f) \|_\infty
\langle \bGa_t (\bP_s \wh{g} - \wh{g}, \bP_s \wh{g} - \wh{g}) \rangle_\mus \,.
\label{eq:gbd1}
\eea 
As $(f,f) \in \DD(\bGa)$ we have the convergence of 
$\bGa_t(f,f)$ to $\bGa(f,f)$ in  
supremum norm. Further, the $\mus$-invariance of $\bP_s$ 
results with 
\begin{equation}\label{eq:gbd2}
\langle \bGa_t (\bP_s \wh{g} - \wh{g}, \bP_s \wh{g} - \wh{g}) \rangle_\mus
\to 
- 2 \langle (\bP_s \wh{g} - \wh{g}) \LL (\bP_s \wh{g} - \wh{g}) \rangle_{\mus}
\end{equation}
as $t \to 0$, which as we have already shown 
converges to zero when $s \to 0$.
Combining all these facts we get the
stated strong continuity of $\bA \bP_s g$.

Fixing $\wh{g}=\bP_v g$ for
$v \geq 0$ and $g \in \GG$, it thus remains only to show that    
for $\bP_t^\delta := (\oo{\bP}_1^{\delta f})_t$ of Proposition
\ref{prop2c},
$$
\rho^\de_t := \de^{-1} (\bP_t^\de - \bP_t) \wh{g}
- \int_0^t \bP_{t-u} \bA \bP_{u} \wh{g} du \,, 
$$
is such that $t^{-1} \| \rho^\de_t \| \to 0$ as $\de \downarrow 0$ 
uniformly in $t \in (0,T]$. 
To this end, recall that the generator $\LL^\de := \oo{\LL}^{\de f}_1$
of the semi-group $\bP_t^\de$ of Proposition \ref{prop2c} has exactly 
the same domain as $\oo{\LL}$, so with
$\wh{\GG} \subseteq \DD(\oo{\LL})$, we see that
both $\bP_t \wh{g}$ and $\bP_t^{\de} \wh{g}$
are in the domain of $\oo{\LL}$ (and $\LL^\de$), with   
$$
\psi_t^\de = \de^{-1} (\LL^\de \bP_t^{\de} \wh{g}
- \oo{\LL} \bP_t^{\de} \wh{g})
- \bA \bP_t \wh{g} \,,
$$
well defined. Further, as we shall prove at the sequel, 
\begin{lemma}\label{lem3}
For any $f \in \GG$, $\wh{g} \in \wh{\GG}$ and $T<\infty$ there exists
$\wh{K}=\wh{K}(f,\wh{g},T)$ finite, such that 
$\|\psi_t^\de\| \leq \wh{K} \sqrt{\de}$ for all 
$t \in [0,T]$ and $\de \in (0,1]$.
\end{lemma}
Next, given Lemma \ref{lem3}, recall that 
$$
\bP_t^\de \wh{g} - \bP_t \wh{g} = \int_0^t \LL^\de \bP_u^\de \wh{g} du -
\int_0^t \oo{\bP}_{t-u} \oo{\LL} \wh{g} du \,,
$$
from which we conclude by the finiteness of $\int_0^t \|\bA \bP_u \wh{g}\| du$
and $\int_0^t \|\psi_u^\de \| du$, that 
$$
\Delta := \de [ \int_0^t  \oo{\bP}_{t-u} \psi_u^\de du - \rho_t^\de ] =
\int_0^t (\oo{\bP}_{t-u}-\bI) \LL^\de \bP_u^{\de} \wh{g} du 
- \int_0^t \oo{\bP}_{t-u} \oo{\LL} (\bP_u^{\de} \wh{g}-\wh{g}) du \,.
$$ 
Fixing $t \in [0,T]$ and $\delta>0$ we claim that $\Delta=
{\bf 0}$. 

Indeed, fixing $h \in \GG$, by Fubini's theorem and the definition of
the adjoint semi-group $(\oo{\bP}^\star)_{t \geq 0}$, clearly, 
$$
\langle h \Delta \rangle_{\mus} = 
\int_0^t \langle (\oo{\bP}^\star_{t-u} h - h) \LL^\de \bP_u^{\de} \wh{g} 
\rangle_{\mus} du 
-
\int_0^t \langle (\oo{\LL}^\star \oo{\bP}_v^\star h) 
(\bP_{t-v}^{\de}\wh{g}-\wh{g}) \rangle_{\mus} dv 
\,.
$$
Further, with $h \in \DD(\oo{\LL}^\star)$ and 
$\wh{g} \in \DD(\LL^\de)$, we thus get that also 
$$
\langle h \Delta \rangle_{\mus} = 
\int_0^t \langle 
(\int_0^{t-u} \oo{\bP}^\star_v \oo{\LL}^\star h dv) 
\LL^\de \bP_u^{\de} \wh{g} \rangle_{\mus} du 
-
\int_0^t \langle \oo{\bP}^\star_v \oo{\LL}^\star h 
(\int_0^{t-v} \LL^\de \bP_u^{\de} \wh{g} du) \rangle_{\mus} dv 
\,.
$$
By Fubini's theorem, the contractiveness of 
$\oo{\bP}^\star_v$ and $\bP^\de_u$ and the finiteness of 
$\| \oo{\LL}^\star h \|$ and $\| \LL^\de \wh{g} \|$, this 
implies that $\langle h \Delta \rangle_{\mus} =0$ for all
$h \in \GG$. So, with $\GG$ dense in $L_2(\mus)$, 
we deduce that $\Delta=
{\bf 0}$ as claimed, i.e.  
$\rho_t^\de = \int_0^t \oo{\bP}_{t-u} \psi_u^\de du$. 

Now, by Lemma \ref{lem3}, 
the contractiveness of $\oo{\bP}_t$ and the convexity of 
$\| \cdot \|$ we have that 
$$
\|\rho_t^\de \| \leq \int_0^t \| \psi_u^\de \| du  \leq t \wh{K} \sqrt{\de}
\,,
$$
thus completing the proof of the theorem.
\qed

\noindent{\bf Proof of Lemma \ref{lem3}:} Fixing 
$f, h \in \GG$, for each
$r \in \DD(\oo{\LL}) \subseteq \DD(\wt{\LL})$ set 
$\xi(r) = \frac{1}{2} \xi_1(r) + \xi_2(r)$, where 
\begin{equation}\label{eq:xidec}  
\xi_1 (r) = 
\langle f r \oo{\LL}^\star h - h f \oo{\LL} r - h r \oo{\LL} f 
\rangle_\mus \quad \mbox{and} \quad
\xi_2 (r) = \langle h f (\wt{\LL} r  - \oo{\LL} r) \rangle_\mus \,. 
\end{equation}
Next for each $r \in \DD(\oo{\LL})$ and $\de \in (0,1]$ let  
$$
\Delta_{\de} (r) 
:= \de^{-1} \langle h (\LL^\de r -\oo{\LL} r) \rangle_\mus - \xi(r) \,.
$$
It is not hard to verify that   
$\xi(r) = \langle h \bA r \rangle_\mus$ whenever   
$r \in \DD(\bA)$. In particular, 
with $\bP_t \wh{g} \in \DD(\bA)$ and $\bP_t^\de \wh{g} \in \DD(\oo{\LL})$ 
this implies that 
\begin{equation}\label{eq:rt-def}
\langle h \psi^\de_t \rangle_\mus = \Delta_\de (\bP_t^\de \wh{g}) 
+ \xi(\bP_t^\de \wh{g}- \bP_t \wh{g}) \,,
\end{equation} 
for any $\wh{g} \in \wh{\GG}$, $\de \in (0,1]$ and $t \in [0,T]$.

To complete the proof we require the following lemmas, whose proofs
are provided at the end of the section.
\begin{lemma}\label{bds1-lem}
Under the conditions of Theorem \ref{thm3b} 
there exists $K_1=K_1(f)<\infty$ such that
for all $\de \in (0,1]$, $h \in \GG$ and $r \in \DD(\oo{\LL})$
\begin{eqnarray}\label{bd-xi}
|\xi(r)| &\leq& K_1 \oo{\EE} (r)^{1/2} \|h \| \,, \\
\label{bd-del}
|\Delta_\de(r)| &\leq& \de K_1 \oo{\EE} (r)^{1/2} \|h\|  \,,
\end{eqnarray}
where $\oo{\EE}(r):=- \langle r \oo{\LL} r \rangle_{\mus}$ is finite 
and non-negative for $r \in \DD(\oo{\LL})$.
\end{lemma}
\begin{lemma}\label{bds-lem}
Under the conditions of Theorem \ref{thm3b}, 
for each $T<\infty$ and $\wh{g} \in \wh{\GG}$ 
there exists $\kappa=\kappa(f,\wh{g},T)<\infty$
such that for all $\de \in (0,1]$ and any $t \in [0,T]$,
\begin{eqnarray}\label{bd-EP}
\oo{\EE} (\bP_t^{\de} \wh{g}) &\leq& \kappa \,, \\
\label{bd-EDP}
\oo{\EE} (\bP_t^{\de} \wh{g}- \bP_t \wh{g}) &\leq& \de \kappa \,.
\end{eqnarray} 
\end{lemma}
Indeed, in view of (\ref{eq:rt-def}) we get from the bounds 
of Lemma \ref{bds1-lem} that
$$
|\langle h \psi^\de_t \rangle_\mus| \leq K_1 \| h \| \big[ \,
\de \oo{\EE} (\bP_t^\de \wh{g})^{1/2} + 
\oo{\EE} (\bP_t^\de \wh{g}- \bP_t \wh{g})^{1/2} 
\big] \,.
$$
Hence, by the bounds of Lemma \ref{bds-lem} we deduce that 
$|\langle h \psi_t^\de \rangle_\mus| \leq \wh{K} \| h \| \sqrt{\de}$ 
for $\wh{K}=\wh{K}(f,\wh{g},T)= 2 K_1 \sqrt{\kappa}$ finite, 
and all $\de \in (0,1]$, $t \in [0,T]$ and $h \in \GG$.
Since $\GG$ is a dense linear subspace of $\HH$
we conclude that $\|\psi_t^\de\| \leq \wh{K} \sqrt{\de}$, as claimed. 
\qed

\noindent{\bf Proof of Lemma \ref{bds1-lem}:} Fixing $f,h \in \GG$ 
and $r \in \DD(\oo{\LL})$ we start by proving (\ref{bd-xi}). 
Indeed, with $\oo{\LL}$ and $\oo{\LL}^\star$
being the generators of the strongly continuous
semi-group $\oo{\bP}_t$ and its adjoint $\oo{\bP}_t^\star$, we have 
from (\ref{eq:xidec}) that
\begin{equation}\label{eq:x1val}
\xi_1(r) = \lim_{t \to 0} \frac{1}{t}  
[ \langle r f \oo{\bP}_t^\star h - h f \oo{\bP}_t r - h r \oo{\bP}_t f + h r f  
\rangle_\mus ]  
= \lim_{t \to 0} \langle h \oo{\bGa}_t (f,r) \rangle_\mus \,,
\end{equation}
for the bi-linear symmetric, non-negative definite operators 
$$
\oo{\bGa}_t(h_1,h_2) 
:=t^{-1} [h_1 h_2 - h_1 \oo{\bP}_t 
h_2 - h_1 \oo{\bP}_t h_2 + \oo{\bP}_t h_1 h_2 ] 
$$
on $\HH \times \HH$. Next, mimicking the arguments of 
(\ref{eq:gbd1}) and (\ref{eq:gbd2}) it follows 
by the $\mus$-invariance of $\oo{\bP}_t$ that 
\begin{equation}
\limsup_{t \to 0} \| \oo{\bGa}_t(f,r) \|^2 
\leq \limsup_{t \to 0} \| \oo{\bGa}_t(f,f) \|_\infty
\langle \oo{\bGa}_t (r,r) \rangle_\mus 
= 2 \| \bGa(f,f) \|_\infty \oo{\EE}(r) \,,
\label{eq:gbd3}
\end{equation} 
for all $r \in \DD(\oo{\LL})$ and $f \in \CCb$ such that 
$(f,f) \in \DD(\bGa)$. Consequently, 
\begin{equation}
|\xi_1(r)| \leq \|h\| \limsup_{t \to 0} \| \oo{\bGa}_t (f,r) \| 
\leq \sqrt{2 \|\bGa (f,f)\|_\infty}\, \oo{\EE}(r)^{1/2} \|h \| \,.
\label{eq:gbd4}
\end{equation} 
Clearly, $|\xi_2(r)| \leq \| h f \|\, \|\wt{\LL} r - \oo{\LL} r \| 
\leq \sqrt{K} \|f \|_\infty \| \oo{\EE}(r)^{1/2} \|h\|$ 
by (\ref{eq:xidec}) and (\ref{eq:sector}). We thus deduce that 
(\ref{bd-xi}) holds for any 
$K_1(f) \geq \sqrt{\|\bGa (f,f)\|_\infty} + \sqrt{K} \|f \|_\infty$.

Turning to (\ref{bd-del}) and fixing $\de \in (0,1]$  
let $\psi_\de := \de^{-1} (e^{\de f} - \de f -\won) \in \CCb$,  
and note that by (\ref{eq:Ladf}) and (\ref{eq:xidec}) 
\begin{equation}\label{eq:dedec}
\Delta_\de(r) = \wh{\Delta}_\de(r) +
\langle h \psi_{-\de} (\wt{\LL} r - \oo{\LL} r) \rangle_\mus
\end{equation}
where 
$$
\wh{\Delta}_\de(r) :=
\de^{-1} [ \langle h \wt{\LL}^{\de f} r \rangle_\mus  
-\langle (1-e^{-\de f}) h \wt{\LL} r \rangle_\mus ]
- \frac{1}{2} \xi_1(r) \,.
$$
Recall that $e^{-\de f} h \in \GG$, so using the 
$\mu^{\de f}$-symmetry of $\wt{\LL}^{\de f}$
and the $\mus$-symmetry of $\wt{\LL}$, followed by the relations  
(\ref{eq:wLdf}) and (\ref{eq:wLd}), 
at $g=e^{-\de f} h$ and $g=h$, we have that 
\beaa
\de [ \wh{\Delta}_\de(r) + \frac{1}{2} \xi_1(r) ] &=&
\langle r e^{\de f} \wt{\LL}^{\de f} (e^{-\de f} h) \rangle_\mus 
-\langle r [\wt{\LL} h- \wt{\LL} (e^{-\de f} h) ] \rangle_\mus \\
&=& \langle r e^{\de f} \wt{\LL} (e^{-\de f} h) \rangle_\mus 
+ \frac{1}{2} \langle r \oo{\bGa}(e^{\de f},e^{-\de f} h) \rangle_\mus 
-\langle r \wt{\LL} h \rangle_\mus \\
&=&
\frac{1}{2} [ \langle r e^{\de f} \oo{\LL}^\star (e^{-\de f} h) \rangle_\mus 
- \langle r h e^{-\de f} \oo{\LL} e^{\de f} \rangle_\mus 
- \langle h \oo{\LL} r \rangle_\mus ] \,.
\eeaa 
Adapting the derivation of (\ref{eq:x1val}) we find that the latter
expression is the limit as $t \to 0$ of 
$\frac{1}{2} \langle h e^{-\de f} \oo{\bGa}_t (e^{\de f},r) \rangle_\mus$,
so in view of (\ref{eq:x1val}), we deduce that 
$$
\wh{\Delta}_\de(r) = \frac{1}{2} \lim_{t \to 0} 
\langle h [ \de^{-1} e^{-\de f} \oo{\bGa}_t (e^{\de f},r) 
- \oo{\bGa}_t (f,r) ] \rangle_\mus \,.
$$
Since $u \mapsto \oo{\bGa}_t(u,r)$ is a linear functional 
such that $\oo{\bGa}_t (\won, r)=0$ for any $r \in \HH$, we further 
have that 
$$
\de^{-1} e^{-\de f} \oo{\bGa}_t (e^{\de f},r) - \oo{\bGa}_t (f,r) 
= e^{-\de f} \oo{\bGa}_t (\psi_\de,r) - \phi_\de \oo{\bGa}_t (f,r)
$$
for $\phi_\de := \won - e^{-\de f}$. With $e^{-\de f} \in [0,1]$ and 
$\psi_\de \in \GG$, we thus get as in the derivation 
of (\ref{eq:gbd3}) and (\ref{eq:gbd4}) that 
\beaa
\| \wh{\Delta}_\de(r) \| &\leq& \frac{1}{2} \| h \| [ \limsup_{t \to 0} 
\|\oo{\bGa}_t (\psi_\de,r)\| 
+ \|\phi_\de\|_\infty \limsup_{t \to 0} \| \oo{\bGa}_t (f,r) \| ] \\
&\leq& k_0(\de,f) \oo{\EE}(r)^{1/2} \| h\| \,,
\eeaa
for $k_0(\de,f)=\sqrt{\|\bGa(\psi_\de,\psi_\de)\|_\infty}
+ \|\phi_\de\|_\infty \sqrt{\| \bGa(f,f) \|_\infty}$ finite.   
Combining this with (\ref{eq:dedec}) and (\ref{eq:sector}) we thus deduce
that 
$$
\| \Delta_\de(r) \| \leq \| \wh{\Delta}_\de(r) \| +
\|h\| \|\psi_{-\de}\|_\infty \|\wt{\LL} r - \oo{\LL} r \| 
\leq k_1(\de,f) \oo{\EE}(r)^{1/2} \| h\| \,,
$$
where $k_1(\de,f) = k_0(\de,f)+ \sqrt{K} \|\psi_{-\de}\|_\infty$.
Thus, we establish (\ref{bd-del}) once we show that
$k_1(\de,f) \leq \de K_1(f)$ for some finite $K_1(f)$ and
all $\de \in (0,1]$. To this end,  
recall that $\| \phi_\de\|_\infty \le \de \|f\|_\infty$ and
$\| \psi_{-\de} \|_\infty \leq \de \|f\|_\infty^2 /2$ for 
any non-negative $f \in \CCb$ and $\de>0$.  
Further, for any $g \in \GG$, as $t \to 0$
the non-negative functions 
$\bGa_t(g,g)$ converge to $\bGa(g,g)$ with respect to the 
supremum norm, so it remains only to check that 
for any $\de \in (0,1]$ and $t>0$,
\begin{equation}\label{eq:gtpsibd}
\bGa_t (\psi_\de,\psi_\de) \leq (e^{\de \|f\|_\infty}-1)^2 \bGa_t (f,f)
\end{equation}
(indeed, using the bound $e^{\de u} - 1 \leq \de e^u$ in (\ref{eq:gtpsibd})
results with (\ref{bd-del}) holding when 
$K_1(f) \geq (e^{\|f\|_\infty} + \|f\|_\infty) 
\sqrt{\|\bGa (f,f)\|_\infty} + \sqrt{K} \|f \|_\infty^2/2$).
As for (\ref{eq:gtpsibd}), recall that for any $x,y \in \sss$ and
$\de>0$, 
$$
| \psi_\de(y)-\psi_\de(x) | = | \int_{f(x)}^{f(y)} (e^{\de u}-1) du | 
\leq c_\de |f(y)-f(x)| \,,
$$
for $c_\de = e^{\de \|f\|_\infty}-1$. Consequently, 
$$
\bGa_t(\psi_\de,\psi_\de) = t^{-1} \E_x
[ ( \psi_\de (X(t)) - \psi_\de(x) )^2 ] \leq 
c_\de^{2} t^{-1} \E_x [ (f(X(t))-f(x))^2 ] =
c_\de^2 \bGa_t(f,f) \,,
$$
which is exactly the inequality (\ref{eq:gtpsibd}).
\qed

\noindent{\bf Proof of Lemma \ref{bds-lem}:} Fixing $T<\infty$,
$f \in \GG$ and $\wh{g} \in \wh{\GG}$, 
recall that $(h,\wh{g}) \in \DD(\bGa) \subseteq \DD(\oo{\bGa})$ 
for any $h \in \GG$. Consequently, then 
$\oo{\bGa}_t(h,\wh{g}) \to \bGa(h,\wh{g})$ 
as $t \to 0$.
In view of (\ref{eq:gbd3}) and (\ref{eq:gtpsibd}) we see that 
$$
\| \oo{\bGa} (\psi_\de,\wh{g}) \|^2  \leq 
2 \| \bGa(\psi_\de,\psi_\de) \|_\infty \, \oo{\EE}(\wh{g})
\leq 2 (e^{\de \|f\|_\infty}-1)^2 \|\bGa(f,f)\|_\infty \, \oo{\EE}(\wh{g})
$$
are  bounded uniformly in $\de \in (0,1]$. Further, since
$e^{\de f} = \de \psi_\de + \de f + \won$, 
\begin{equation}\label{eq:unif-bd}
\| \oo{\bGa} (e^{\de f},\wh{g}) \|
\leq \de \| \oo{\bGa} (\psi_\de,\wh{g}) \| + \de \| \oo{\bGa} (f, \wh{g}) \|
\leq \kappa_0(f) \oo{\EE}(\wh{g})^{1/2}
\end{equation} 
are also bounded uniformly in $\de \in (0,1]$ (take  
$\kappa_0(f) = e^{\|f\|_\infty} \sqrt{2 \|\bGa(f,f)\|_\infty}$ finite).

Turning to prove (\ref{bd-EP}), note that for any 
$\de \in (0,1]$ and $r \in \DD(\oo{\LL})$, by (\ref{eq:Ladf}),
$$
-\langle r \oo{\LL}^{\de f}_1 r \rangle_{\mu^{\de f}}
= 
- \langle r \wt{\LL}^{\de f} r \rangle_{\mu^{\de f}}
- \langle r \oo{\LL} r \rangle_{\mus} \geq \oo{\EE} (r)  
$$
since 
$\wt{\LL}^{\de f}$ is the generator of a 
$\mu^{\de f}$-symmetric, strongly continuous semi-group
whose domain contains $\DD(\oo{\LL})$, hence a negative
self-adjoint operator on this set. 
With $\bP^\de_t$ a contraction 
for the norm of $L_2(\mu^{\de f})$ (denoted hereafter by 
$\|\cdot\|_\de$), we thus find that 
for any $\wh{g} \in \wh{\GG}$ and $t>0$,
$$ 
\oo{\EE} (\bP_t^{\de} \wh{g}) \leq
-\langle (\bP^\de_t \wh{g})  
(\bP^\de_t \oo{\LL}^{\de f}_1 \wh{g}) \rangle_{\mu^{\de f}}
\leq \| \wh{g} \|_\de  \, \| \oo{\LL}^{\de f}_1 \wh{g} \|_\de \,.
$$ 
Further, $f$ is non-negative and
by our assumptions, $\wh{g}$ is in $\CDD$ of (\ref{eq:whddef}) 
where the identity (\ref{eq:wLdf}) applies.
We thus establish the uniform bound of (\ref{bd-EP}) upon 
noting that both 
$\| \wh{g} \|_\de \le e^{\|f\|_\infty} \|\wh{g}\|$ and
\bea
\| \oo{\LL}^{\de f}_1 \wh{g} \|_\de
&\leq& \| \wt{\LL}^{\de f} \wh{g} \|_\de
+ \| e^{-\de f} \oo{\LL} \wh{g} \|_\de \nonumber \\
&\leq&
\frac{e^{\|f\|_\infty}}{2} (\|\oo{\LL} \wh{g} \| 
             +\| \oo{\LL}^\star \wh{g} \| 
+ \| \oo{\bGa} (e^{\de f},\wh{g}) \| ) 
+ e^{\|f\|_\infty} \| \oo{\LL} \wh{g} \|
\label{eq:bdLd}
\eea
are bounded uniformly in $\de \in (0,1]$ (see 
(\ref{eq:unif-bd}) for the uniform bound on $\| \bGa (e^{\de f},\wh{g}) \|$).

Next consider the non-negative quadratic form 
$\wt{\EE}(h_1,h_2):= - \langle h_1 \wt{\LL} h_2 \rangle_{\mus}$
on $\DD(\wt{\LL}) \times \DD(\wt{\LL})$, noting that 
$\oo{\EE} (r) = \wt{\EE}(r,r)$ 
for all $r \in \DD(\oo{\LL}) \cap \DD(\oo{\LL}^\star)$,
due to (\ref{eq:wLd}). Turning to prove (\ref{bd-EDP}), 
by our assumptions this is the case for 
$r_t^\de := \bP_t^{\de} \wh{g}- \bP_t \wh{g}$, 
hence $w_t^\de := \oo{\EE}(r_t^\de)=\wt{\EE}(r_t^\de,r_t^\de)$. 
Further, 
$t \mapsto a_t^\de := 
\partial_t r_t^\de = \LL^\de \bP_t^\de \wh{g} - \oo{\LL} \bP_t \wh{g}$ 
is uniformly continuous, since 
\begin{equation}\label{eq:unifc}
\sup_{|t-s| \leq \theta} \,
\|a_t^\de-a_s^\de\| \leq \sup_{u \leq \theta} 
\| (\bP_u^\de - \bI) \LL^\de \wh{g} \| + \sup_{u \leq \theta} 
\| (\bP_u - \bI) \LL \wh{g} \| =: \eps_\de(\theta)
\end{equation}
which converges to zero when $\theta \downarrow 0$. 
With $b_t^\de = a_t^\de - \oo{\LL} r_t^\de$,
we show next that  
\begin{eqnarray}\label{eq:bdudel}
\|b_t^\de\| &\leq& 2 \de K_1 \sqrt{\kappa} \,,\\
\| \oo{\LL} r_t^\de \| &\leq& \kappa_1 \,,
\label{eq:bdlrdel}
\end{eqnarray}
for some finite $\kappa_1=\kappa_1(f,\wh{g})$, the universal constants
$K_1$ and $\kappa$ of 
Lemma \ref{bds1-lem} and (\ref{bd-EP}), respectively, all $t \geq 0$
and $\de \in (0,1]$. Indeed,
$b_t^\de = \LL^\de \bP_t^{\de} \wh{g} - \oo{\LL} \bP_t^\de \wh{g}$,
hence $\langle h b_t^\de \rangle_{\mus} = 
\delta [ \Delta_\de(\bP_t^\de \wh{g}) + \xi(\bP_t^\de \wh{g}) ]$ 
for any $h \in \GG$. Thus, by Lemma \ref{bds1-lem} and (\ref{bd-EP})  
$$
|\langle h b_t^\de \rangle_{\mus}| \le 
\delta K_1 (\de + 1) \|h\| \oo{\EE}(\bP_t^\de \wh{g})^{1/2} 
\leq  \delta K_1 (\de + 1) \sqrt{\kappa} \|h\|  \,,
$$
and with $\GG$ dense in $L_2(\mus)$ this immediately yields the 
bound of (\ref{eq:bdudel}). Turning to prove (\ref{eq:bdlrdel}),
note that for any $t \geq 0$ and $\de \in (0,1]$, by (\ref{eq:bdudel})
the contractivity of $\bP_t$ on $L_2(\mus)$ and the contractivity of $
\bP_t^\de$ on $L_2(\mu^{\de f})$,
\beaa
\| \oo{\LL} r_t^\de \| &\leq& \| \LL^\de 
\bP_t^\de \wh{g} \| + \| b_t^\de \| + \| \oo{\LL} \bP_t \wh{g} \| \\
&\leq& e^{\|f\|_\infty}
\| \LL^\de \wh{g} \|_{\de} + 2 K_1 \sqrt{\kappa} + \| \oo{\LL} \wh{g} \| \,,
\eeaa 
with the right side bounded uniformly in $\de \in (0,1]$ by
some finite $\kappa_1=\kappa_1(f,\wh{g}) \ge 1$ 
(see (\ref{eq:bdLd}) for details).

As $\wt{\LL}$ is $\mus$-symmetric, by Fubini's theorem, for  
any $s' \geq s$,
$$
w_{s'}^\de-w_s^\de -(s'-s) 
[\wt{\EE}(a_s^\de,r_s^\de) + \wt{\EE}(a_{s'}^\de,r_{s'}^\de)]= 
\int_s^{s'} \wt{\EE}(a_u^\de-a_s^\de,r_s^\de) du 
- \int_s^{s'} \wt{\EE}(a_{s'}^\de - a_u^\de,r_{s'}^\de) du \,. 
$$
Recall that by (\ref{eq:sector}) and (\ref{eq:bdlrdel}),
for some universal finite constant $K$,
\begin{equation}\label{eq:wtllbd}
\| \wt{\LL} r_t^\de \| \leq \sqrt{K w_t^\de} + \|\oo{\LL} r_t^\de\|
\leq K w_t^\de + 2 \kappa_1 \,,
\end{equation}
so by the uniform continuity of $t \mapsto a_t^\de$, see (\ref{eq:unifc}),
$$
w_{s'}^\de-w_s^\de -(s'-s) 
[\wt{\EE}(a_s^\de,r_s^\de) + \wt{\EE}(a_{s'}^\de,r_{s'}^\de)] \le 
(s'-s) \eps_\de (s'-s) [K w_{s'}^\de + K w_s^\de + 4 \kappa_1] \,.
$$
Further, by (\ref{eq:sector}) also 
$$
\wt{\EE}(a_t^\de-b_t^\de,r_t^\de) =
- 2 \langle (\oo{\LL} r_t^\de) \wt{\LL} r_t^\de \rangle_{\mus} \le 
\|\wt{\LL} r_t^\de  - \oo{\LL} r_t^\de \|^2 
\leq K \oo{\EE} (r_t^\de) = K w_t^\de \,.
$$
Combining the latter pair of inequalities and the bound 
$\wt{\EE}(b_t^\de,r_t^\de) \leq \|b_t^\de\| \| \wt{\LL} r_t^\de\|$,
we deduce from (\ref{eq:bdudel}) and (\ref{eq:wtllbd})
that if $\eps_\de(s'-s) \leq 2 \de K_1 \sqrt{\kappa}$, then 
$$
w_{s'}^\de-w_s^\de \le  
(s'-s) [K_2 w_{s'}^\de+ K_2 w_s^\de + \kappa_2 \de ] \,,
$$
for some finite universal constants
$K_2=K_2(K_1,\kappa,K)$, $\kappa_2=\kappa_2(K_1,\kappa,\kappa_1)$
and all $\de \in (0,1]$.
Since $w_0^\de=0$, iterating the latter bound $n$ times, for 
$s'-s=t/n$, then taking $n \to \infty$, we conclude that 
$w_t^\de \leq \de \kappa_2 \int_0^t e^{2K_2 u} du$ for all 
$\de \in (0,1]$ and $t \geq 0$. That is, 
(\ref{bd-EDP}) holds for $\kappa(f,\wh{g},T)=\kappa_2 \int_0^T e^{2K_2 u} du$ 
finite.
\qed

\section{Pure jump processes on a discrete state space}\label{sec:jump}

We consider here pure 
jump processes $X(t)$ on a countable (or finite) 
state space $\sss$ equipped with the discrete topology,
such that the total jump rate at state $x$ is bounded 
uniformly over $x \in \sss$. That is, the 
jump rates $c(x,y)\ge 0$ from $x$ to $y \neq x$ 
are such that 
\begin{equation}\label{eq:jd1}
\sup_{x\in \sss} \sum_{y:y\ne x} c(x,y)<\infty
\end{equation}
(which trivially holds when the set $\sss$ is finite).
Recall that to each such process corresponds a strongly continuous Markov
semi-group on the Banach space
$\CC_b$ of all bounded functions on $\sss$,
the generator of which is the bounded linear
operator $\LL: \CC_b \to \CC_b$ such that  
\begin{equation}\label{eq:jd0}
\LL g(x)=\sum_{y:y \neq x} c(x,y)(g(y)-g(x)) \,.
\end{equation}
Conversely, any operator of the form (\ref{eq:jd0})
with non-negative $c(x,y)$ satisfying (\ref{eq:jd1})
is the generator of such a Markov process, and 
taking in this context $\GG=\CC_b$ eliminates
all technical issues of the previous sections
(about the domains of various generators).
Further assuming that the process 
$X(\cdot)$ is irreducible, let $\mus(\cdot)$
denote its unique invariant probability measure, identified
hereafter with the positive function $\mus(x):=\mus(\{x\})$ on $\sss$.
Recall that necessarily, 
\begin{equation}\label{eq:jd2}
\mus(x)\sum_{y:y\ne x} c(x,y)=\sum_{y:y\ne x} \mus(y)c(y,x)\,,
\qquad \forall x\in \sss \,.
\end{equation}

We proceed to compute in this case 
the response functions for our two 
generic Markov perturbations. 
To this end, consider the operators $\LL^\star$ and 
$\wt{\LL} = \frac{1}{2} (\LL + \LL^\star)$ of the
form (\ref{eq:jd0}) but for jump rates 
\beaa
c^\star(x,y)&:=&\frac{\mus(y)}{\mus(x)}c(y,x),\qquad x\ne y\,,\\
\wt{c}(x,y)&:=&\frac{1}{2}(c(x,y)+c^\star(x,y))
\,\qquad x \ne y \,,
\eeaa
respectively. By (\ref{eq:jd2}), both $c^\star(\cdot,\cdot)$ and 
$\wt{c}(\cdot,\cdot)$ satisfy (\ref{eq:jd1}) so 
$\LL^\star$ and $\wt{\LL}$ are both bounded operators on $\CC_b$
and the generators of strongly continuous, 
Markov semi-groups on $\CC_b$, denoted $\bP_t^\star$ and $\wt{\bP}_t$,
respectively. It is easy to check that $\bP_t^\star$ 
is the adjoint of the original semi-group $\bP_t$
and that $\wt{\bP}_t$ is  
the $\mus$-symmetric, strongly continuous, Markov semi-group 
of Corollary \ref{cor2b} (both restricted
to $\CC_b \subseteq L_2(\mus)$). 

Given a non-negative $f\in\CC_b$,
both Propositions \ref{prop1} and \ref{prop2c} 
apply here, and their generic perturbations correspond  
to the bounded operators $\LL^f_0$ and $\LL^f_1$ on $\CC_b$ 
having jump rates $c^{f}_0(x,y)=e^{-f(x)}c(x,y)$ and
\beaa 
c^f_1(x,y)&=& \wt{c}(x,y)+
\frac{1}{2}(e^{f(y)-f(x)}-1)c(x,y)+\frac{1}{2}e^{-f(x)}(c(x,y)-c^\star(x,y))\\
&=&
\frac{1}{2}(e^{f(y)-f(x)}+e^{-f(x)})c(x,y)+\frac{1}{2}(1-e^{-f(x)})c^\star(x,y),
\eeaa 
respectively. Theorem \ref{thm3a} provides the 
response function for $X_0^f(\cdot)$ which in this case
has a bounded operator $(\bA_0)_f$ on $\CC_b$ 
of the form
\begin{equation}\label{eq:jd3}
(\bA_f) g (x)=\sum_{y:y\ne x} a^f(x,y) (g(y)-g(x)) \,,
\end{equation}
with $a^f_0(x,y)=-f(x)c(x,y)$. 
Consider the bounded operator $(\bA_1)_f$ of the form (\ref{eq:jd3}) 
with 
$$
a^f_1(x,y)=
a^f_0(x,y)+\frac{1}{2}f(y)c(x,y)+\frac{1}{2} f(x) c^\star(x,y) \,,
$$
and let
\beaa
\xi_\de(x,y)&:=& \de^{-1} (c_1^{\de f}(x,y)-c(x,y)) - a^f_1(x,y) \\
&=& \frac{1}{2} \varphi_{-\de}(f(x)) [c^\star(x,y)-c(x,y)] +
\frac{1}{2} \varphi_\de(f(y)-f(x)) c(x,y) \,,
\eeaa
where $\varphi_\de(r) = \de^{-1} (e^{\de r}-\de r -1) \to 0$
as $\de \to 0$, uniformly on compacts. Hence,
\begin{equation}\label{eq:xibd}
\lim_{\de \to 0} \sup_{x \in \sss} \sum_{y: y \neq x} |\xi_\de(x,y)| = 0 \,.
\end{equation}
This in turn implies that (\ref{eq:am1}) holds for $\LL^f_1$ and 
$(\bA_1)_f$, so by Proposition \ref{prop0} we deduce that 
the Markovian perturbations $X_1^f(\cdot)$ have the response
function associated with $(\bA_1)_f$.

Suppose now that $\bA_f$ is of the form (\ref{eq:jd3}) and that
\begin{equation}\label{eq:jd4a}
\sup_{x\in \sss} \sum_{y:y\ne x} |a^f(x,y)|<\infty 
\end{equation}
which guarantees that $\bA_f$ is a bounded operator on $\CC_b$.
In view of Theorem \ref{theo-fdt}, for such $\bA_f$ to correspond
to the response function of some perturbation $X^f(\cdot)$ 
(per Definition \ref{def:resp}), it is necessary that 
$a^f(x,y)=b^f(x,y)-f(x) c(x,y)$, where
\beqn{eq:jd4}
b^{rf}(x,y)&=&rb^f(x,y),\quad r>0 \\
\label{eq:jd5}
\mus(x) \sum_{y:y\ne x} b^f(x,y)&=&\sum_{y:y\ne x} 
\mus(y) b^f(y,x),\qquad \forall x\in \sss \,.
\eeqn 
Slightly modifying the time change 
generic perturbation of Proposition \ref{prop1}, 
we next show that essentially these conditions on $a^f(x,y)$
are also sufficient
for having a perturbation $X^f(\cdot)$ whose response function 
is given by (\ref{res:form}).
\begin{proposition}\label{prop6}
Suppose that the generator of the semi-group $\bP_t$ of a
pure jump Markov process $X(\cdot)$ on a discrete 
state space $\sss$ is of the form (\ref{eq:jd0}) 
for jump rates $c(x,y) \geq 0$ that satisfy  
(\ref{eq:jd1}). To any $\bA_f$ of the form (\ref{eq:jd3}) with 
$a^f(x,y)=b^f(x,y)-f(x) c(x,y)$ satisfying  
(\ref{eq:jd4a})--(\ref{eq:jd5}) and such that 
for some $\rho^f < \infty$, 
\begin{equation}\label{eq:jd5b}
b^f(x,y) \geq -\rho^f c(x,y) \qquad \forall x \ne y \,,
\end{equation}
there corresponds a Markovian perturbation $X^f(\cdot)$ satisfying
Assumption \ref{ass:1} whose response function is 
$\Rfg (s,t)=\bP_s \bA_f \bP_{t-s} g$ (for $g \in \CC_b$).
\end{proposition}
\begin{remark} Condition (\ref{eq:jd5b}) implies that 
$b^f(x,y)$ is non-negative for every $x \neq y$ such 
that $c(x,y)=0$ (and for a finite state space $\sss$
it puts no other restrictions on $b^f(x,y)$). 
\end{remark}

\noindent
\Proof In view of (\ref{eq:jd4}) and 
(\ref{eq:jd5b}), if $\de>0$ is small enough 
so $1-\de \rho^f > 0$ then 
\begin{equation}\label{eq:jd6}
\wh{c}_0^{\delta f}(x,y)
=e^{-\delta f(x)}\Big[c(x,y)+b^{\delta f}(x,y)\Big]=c_0^{\delta f}(x,y)
+e^{-\delta f(x)}b^{\delta f}(x,y),
\end{equation}
are non-negative for all $x \neq y$. Further, by the boundedness of 
$f(\cdot)$, (\ref{eq:jd1}) and (\ref{eq:jd4a}) we have that 
$\sum_{y} \wh{c}_0^{\delta f}(x,y)$ is bounded, uniformly in $x \in \sss$.
Thus, there exists a pure jump Markov 
process $\wh{X}_0^{\de f}(\cdot)$ on $\sss$
whose semi-group is generated by a bounded 
operator $\wh{\LL}_0^{\de f}$ on $\CC_b$ of the form (\ref{eq:jd0})
with the jump rates $\wh{c}_0^{\de f}(x,y)$. Moreover, it follows from
(\ref{eq:jd2}), (\ref{eq:jd5}) and (\ref{eq:jd6}) that 
$$
\mu^{\de f} (x) \sum_{y:y\ne x} \wh{c}_0^{\de f} (x,y)=\sum_{y:y\ne x} 
\mu^{\de f} (y) \wh{c}_0^{\de f} (y,x),\qquad \forall x\in \sss 
$$
(where $\mu^{\de f}(x)=e^{\de f(x)} \mus(x)$). This implies
that $\mu^{\de f}(\cdot)$ is a 
finite, positive invariant measure for the semi-group of
the irreducible Markov process $\wh{X}_0^{\de f}(\cdot)$,
which thus satisfies Assumption \ref{ass:1} (with 
$\GG=\CC_b=\BB$). It is easy to check that (\ref{eq:xibd}) holds 
for $\xi_\de(x,y):= \de^{-1} (\wh{c}_0^{\de f}(x,y)-c(x,y)) - a^f(x,y)$.
This in turn implies that (\ref{eq:am1}) holds in this setting,
so by Proposition \ref{prop0} we deduce that (\ref{res:smpl}) holds
as well, and that $\bA_f \bP_s g$ is strongly continuous on $\CC_b$.
\qed

\begin{remark} We 
alternatively get the response
function $\Rfg (s,t)=\bP_s \bA_f \bP_{t-s} g$ per Proposition \ref{prop6}
by adapting instead the generic perturbation of 
Proposition \ref{prop2c}, i.e.
following the same line of reasoning for 
the Markov perturbation $\wh{X}_1^{\de f}(\cdot)$ that correspond to 
the jump rates 
$$
\wh{c}_1^{\delta f}(x,y)=c_1^{\delta f}(x,y)
+e^{-\delta f(x)}\Big(b^{\de f}(x,y)
-\frac{\de}{2} f(y)c(x,y)-\frac{\de}{2} f(x)c^\star(x,y)\Big) \,.
$$
Indeed, for $f(\cdot)$ non-negative and $\de>0$ also 
$\psi_\de = \de^{-1} (e^{\de f} - \de f -\won)$ is non-negative, so  
by (\ref{eq:jd4}) and (\ref{eq:jd5b})
\beaa
e^{\de f(x)} \wh{c}_1^{\de f}(x,y) &=&
c(x,y)+ b^{\de f} (x,y) + 
\frac{\de}{2} [\psi_\de(y) c(x,y) +\psi_\de(x) c^\star(x,y)] \\
&\ge& c(x,y) + \delta b^f(x,y) 
\eeaa
is non-negative as soon as $1-\de \rho^f > 0$. 
With $\psi_\de \in \CC_b$ it follows from 
(\ref{eq:jd1}) and (\ref{eq:jd4a}) that 
$\wh{\LL}_1^{\de f}$ of the form 
(\ref{eq:jd0}) corresponding to jump rates $\wh{c}_1^{\de f}(x,y)$
is a bounded operator on $\CC_b$ hence a generator of a semi-group
for a Markov process $\wh{X}_1^{\de f}(\cdot)$. Further,  
from the $\mus$-invariance of 
$\LL$ and $\bB_f=\bA_f+f \LL$, see (\ref{eq:jd2}) and (\ref{eq:jd5}),
it follows that 
$\langle \wh{\LL}_1^{\de f} g \rangle_{\mu^{\de f}} =0$
for all $g \in \CC_b$ and so 
the irreducible Markov process $\wh{X}_0^{\de f}(\cdot)$
satisfies Assumption \ref{ass:1}. Finally, with 
$\| \psi_\de\|_\infty \to 0$ as $\de \downarrow 0$
we get the stated response 
function upon checking that (\ref{eq:xibd}) holds for  
$\xi_\de(x,y)= \de^{-1} (\wh{c}_1^{\de f}(x,y)-c(x,y)) - a^f(x,y)$.
\end{remark}

\medskip
Cycle decomposition provides a canonical construction 
of Markov processes on a discrete state space with a
prescribed invariant measure (such as $\mu^f$).
For simplicity, we consider only cycles of
finite length. More precisely, 
equipping $\sss$ with any complete order, 
let $\Gamma$ denote the collection of all  
finite oriented cycles $\gamma$. That is,
$\gamma=(x_0,x_1,\ldots,x_n)$ of length 
$n=|\gamma| \geq 2$ is such that $x_n=x_0$ and $x_i \neq x_j$
for all $0 \le i < j <n$. Suppose 
a strictly positive probability measure $\mus$ on $\sss$ and
$\alpha:\Gamma\mapsto \reals_+$ are such that 
$$
\|\alpha\|_\Gamma := \sup_{x \in \sss} \frac{1}{\mus(x)} 
\sum_{\gamma: x \in \gamma} |\alpha(\gamma)|  
$$ 
is finite (in particular, if $\sss$ is finite then so is 
$\Gamma$ and $\|\alpha\|_{\Gamma} < \infty$ for any
$\alpha : \Gamma \mapsto \reals$).
It is then easy to check that the jump rates
\begin{equation}\label{eq:jd7}
c(x,y)=\frac{1}{\mus(x)}\sum_{\gamma\in\Gamma}
\alpha(\gamma)1_{(x,y)\in\gamma} \qquad \forall x \neq y \in \sss\,,
\end{equation}
satisfy (\ref{eq:jd1}) and that $\mus(\cdot)$ 
is an invariant measure for the corresponding semi-group $\bP_t$.
Further, this semi-group is $\mus$-symmetric (i.e.
the Markov process is reversible), if 
$\alpha(\gamma)=0$ whenever $|\gamma|>2$.
Next, 
let $\beta^f:\Gamma \mapsto \reals$ 
and $\alpha^f: \Gamma \mapsto \reals_+$
for $f\in \CC_b$ be such that $\|\beta^f\|_\Gamma$ is finite,
$\beta^{r f}=r\beta^f$ and 
$$
\lim_{\de \downarrow 0}
\| \de^{-1} (\alpha^{\de f} - \alpha) - \beta^f \|_\Gamma = 0 \,.
$$ 
Then, applying once more Proposition \ref{prop0},
we deduce that the Markov process of jump rates
$$
c^{f}(x,y)=\frac{e^{-f(x)}}{\mus(x)}
\sum_{\gamma\in\Gamma}\alpha^{f}(\gamma)1_{(x,y)\in\gamma}$$
has an invariant measure $\mu^{f}$ and the response
function corresponding to 
$$a^f(x,y)=-f(x)c(x,y)+\frac{1}{\mus(x)}\sum_{\gamma\in\Gamma}
\beta^f(\gamma) 1_{(x,y) \in\gamma} \, .$$

We next consider some concrete examples such as 
the Glauber and Metropolis dynamics for Gibbs 
measures on finite graphs. 
\begin{example}\label{exmp1}
Consider a finite graph with $\sss$ denoting its 
vertexes and the symmetric 
$E \subseteq \sss \times \sss$
denoting its edges. Given $H: \sss \mapsto \reals$ 
let $\mus(x)=e^{-H(x)}$ denote the
corresponding non-normalized Gibbs measure and
consider the reversible Markov processes
obtained by (\ref{eq:jd7}) when  
$\alpha(\gamma)>0$ if and only if 
$\gamma=(x,y,x)$ with $(x,y)\in E$. That is,
having jump rates 
$$
c(x,y)=e^{H(x)} \alpha(x,y) 1_{(x,y) \in E} \,,
$$
for $\alpha(x,y)=\alpha(y,x) > 0$. Two  
such examples are the Metropolis dynamics where  
$\alpha_M (x,y)=\min(e^{-H(x)},e^{-H(y)})$ and 
the Glauber dynamics for which 
$\alpha_G (x,y)=1/(e^{H(x)}+e^{H(y)})$. 
The Markov perturbations one uses for  
the Metropolis (or Glauber) dynamics are 
of the same type as the original process,
just replacing $H(\cdot)$ by $H(\cdot)-\de f(\cdot)$. 
Here the convergence in $\|\cdot\|_\Gamma$ is 
equivalent to a point-wise convergence on $E$
leading to the response functions that correspond to 
\beaa
a_M^f(x,y)&=&(f(y)-f(x)) e^{-\Delta H} \big(
1_{\Delta H>0} + 1_{\Delta H = 0} 1_{f(x)>f(y)} \big) 1_{(x,y)\in E} \,,\\
a_G^f(x,y)&=&(f(y)-f(x)) \frac{e^{-\Delta H}}{\big(1+e^{-\Delta H}\big)^2} 1_{(x,y) \in E} \,,
\eeaa
where $\Delta H := H(y)-H(x)$. Note in particular that while $a_G^f$ is linear in $f(\cdot)$,
this is in general not the case for $a_M^f$.
\end{example}

\section{Finite dimensional diffusion processes}\label{sec:diff}

Here $\sss$ is a connected,
finite dimensional 
$\CC^\infty$-manifold $M$ without boundary. 
We first consider compact $M$, with the treatment of 
non-compact $M=\Bbb R^d$ provided at the end of the section.
Let $\mus$ be a probability 
measure on $M$ that has a smooth strictly positive 
density with respect to any coordinate chart for $M$
(c.f. \cite[Section 6.3]{De-St}).
Setting $\GG=\CC^\infty(M)$ and 
$\Gamma({\bf T}(M))$ denoting the space of smooth sections over $M$, 
recall that in the absence of boundary, for any 
$\bZ \in \Gamma({\bf T}(M))$ there exists then
a unique $g_\bZ \in \GG$ such that 
$\bZ^\star h=-\bZ h + g_\bZ h$ acts on $\GG$ as
the adjoint of $\bZ$ with respect to
the inner product of $L_2(\mus)$ (and 
$g_\bZ = \bZ^\star \won$). 
For $\bX_i \in \Gamma({\bf T}(M))$, $i=0,1,\ldots,d$,
consider the operator   
$\LL=\sum_{i=1}^d \bX_i \circ \bX_i + \bX_0$ on $\GG$ 
(that is, a diffusion generator in the H\"ormander form,
see \cite{hormander}). 
Any such operator can be rewritten as
$$
\LL=-\sum_{i=1}^d \bX^\star_i\circ \bX_i + \bY \,,
$$
where $\bY= \bX_0+\sum_{i=1}^d g_{\bX_i}\,\bX_i$ is 
also in $\Gamma({\bf T}(M))$. Note that such operator is the
restriction to $\GG$ of the generator of a strongly continuous 
Markov semi-group $\bP_t$ on $\CC_b = \CC_b(\sss)$ such that 
$\wh{\GG}=\GG=\CC^\infty(M)$
(see for example, \cite[Theorem 6.3.2]{De-St}). In 
particular, the corresponding Markov process $X(\cdot)$
can be constructed as the unique solution of a certain 
Stratonovich stochastic differential equation
(S-SDE), c.f. \cite[Exercise 6.3.22]{De-St}.

Also, the semi-group $\bP_t$ is $\mus$-invariant if and only if 
$g_\bY = {\bf 0}$ (see \cite[Theorem 6.3.2]{De-St}). In particular,
if $\{\bX_1,\ldots,\bX_d\}$ satisfies H\"ormander's strong hypo-elliptic 
condition (i.e. ({\bf H}) of \cite[Section 6.3]{De-St}), then 
for any $\bX_0 \in \Gamma({\bf T}(M))$
there exists a unique probability measure $\mus$ 
(of a smooth strictly positive 
density with respect to any coordinate chart) for which 
$g_\bY={\bf 0}$. We next provide the diffusion process 
$X^f_1(\cdot)$ of the generalized Langevin dynamics of
Proposition \ref{prop2c} and its response function. 
\begin{proposition}\label{prop:diff}
Suppose $g_\bY={\bf 0}$ and  
for any non-negative $f \in \CC^\infty(M)$ 
let $\bY^f=e^{-f} \bY+\sum_{i=1}^d (\bX_i f) \bX_i$.
Then,
$$
\LL^f_1 =-\sum_{i=1}^d \bX^\star_i\circ \bX_i  + \bY^f \,,
$$
is (the restriction to $\GG=\CC^\infty(M)$ of) the generator of a 
strongly continuous, $\mu^f$-invariant
Markovian semi-group $(\bP^f_1)_t$ on $\CC_b$ with $\GG$
closed under the action of this semi-group  
(and the Markov process $X_1^f(\cdot)$
is the unique solution of a certain
S-SDE). If further 
\begin{equation}\label{eq:njd}
\int_M |\bY g|^2 d\mus \le K \sum_{i=1}^d \int_M |\bX_i g|^2 d\mus \,,
\end{equation}
for some $K<\infty$ and all $g \in \GG$,
then $X^f_1(\cdot)$ has a response function that corresponds to 
$(\bA_1)_f = \sum_{i=1}^d (\bX_i f) \bX_i - f \bY$,
a linear operator of domain $\GG$. 
\end{proposition}
\begin{remark}
The process $X^f_1(\cdot)$ is a diffusion on $M$ that differs  
from $X(\cdot)$ only by the addition of a smooth drift term corresponding to 
$\bY^f-\bY$. We note in passing that $X(\cdot)$ is reversible (i.e. has a 
$\mus$-symmetric semi-group) if and only if
$\bY={\bf 0}$, in which case the added drift
is of a gradient form
(and $X^f_1(\cdot)$ is known in the literature as the Langevin dynamic).
\end{remark}
\begin{remark} If $\bX_0 =\sum_{i=1}^d h_i \bX_i$ for some 
$h_i \in L_2(\mus)$, then 
Girsanov transformation shows that the laws of 
$X^f_1(\cdot)$ and $X(\cdot)$ are mutually absolutely continuous 
on $\CC ([0,T];M)$ for each $T<\infty$. The Langevin dynamic is in this respect
more natural that the time change generic perturbation
(of generator $\LL_0=e^{-f}\LL$), for which this is of course not the case.
If in addition $h_i\in L_\infty(\mus)$ for $i=1,\ldots,d$ 
(as for example, in case of uniform ellipticity), then the
condition (\ref{eq:njd}) is trivially satisfied.
\end{remark}

\noindent \Proof Our assumption that $g_{\bY}={\bf 0}$ means that 
$\bY^\star = -\bY$. Hence, acting on $\GG$, the adjoint 
$$\LL^\star=-\sum_{i=1}^d \bX^\star_i\circ \bX_i - \bY$$
of $\LL$ with respect to the inner product of $L_2(\mus)$ 
is such that $\LL^\star \GG \subseteq \GG$. It is thus 
the restriction to $\GG$ (and a core) 
of the generator of the adjoint semi-group $\bP_t^\star$ 
on $\CC_b$, 
with $\GG$ closed under the action of $\bP_t^\star$. Consequently, 
the generator of the
$\mus$-symmetric semi-group $\wt{\bP}_t$ of Corollary \ref{cor2b} is
just $\wt{\LL}=-\sum_{i=1}^d \bX^\star_i\circ \bX_i$ (when 
acting on $\GG$). Moreover, here 
\begin{equation}\label{eq:gmdiff}
\bGa(f,g)=2\sum_{i=1}^d (\bX_i f) (\bX_i g)
\end{equation} 
satisfies the Leibniz rule $\bGa(fh,g)=h\bGa(f,g)+f \bGa(h,g)$,
so $\bGa(e^f,g)=e^f\bGa(f,g)$ and for any non-negative
$f \in \GG$, the generator of 
the $\mu^f$-symmetric semi-group $\wt{\bP}^f_t$ 
of Lemma \ref{lem2b} is such that for $g \in \GG$,
$$
\wt{\LL}^f g = - (1-e^{-f}) \sum_{i=1}^d \bX^\star_i \circ \bX_i g 
+ \sum_{i=1}^d (\bX_i f) (\bX_i g) \,.
$$
It follows that when acting on $\GG$ the generator 
$\LL^f_1 = \wt{\LL}^f + e^{-f} \LL$ we use in Proposition 
\ref{prop2c} is merely 
$$
\LL^f_1 =-\sum_{i=1}^d \bX^\star_i\circ \bX_i  + \bY^f \,.
$$
It is easy to check that for any $\bZ \in \Gamma({\bf T}(M))$,
the operator $\bZ^{\star,f} = -\bZ + (g_\bZ -\bZ f)$ acts on $\GG$ 
as the adjoint of $\bZ$ with respect to
the inner product of $L_2(\mu^f)$. Further, with
$\bY e^{-f} = -e^{-f} \bY f$, it follows that 
$$
g_{e^{-f} \bY} - e^{-f} \bY f = e^{-f} g_\bY = {\bf 0} \,.
$$
Thus,
$\LL^f_1 = - \sum_{i=1}^d \bX^{\star,f}_i \circ \bX_i + e^{-f} \bY$,
and $(e^{-f} \bY)^{\star,f}=-(e^{-f} \bY)$. Hence, 
by \cite[Theorem 6.3.2]{De-St}, now with respect
to the finite measure $\mu^f$ on $M$ that is also 
of a smooth strictly positive density,
we find that $\LL^f_1$ is the generator of a 
strongly continuous, $\mu^f$-invariant
Markovian semi-group $(\bP^f_1)_t$ on $\CC_b$ such that 
$\GG$ is closed under its action, with $X^f_1(\cdot)$ characterized as
the unique strong solution of 
some S-SDE (this direct construction bypasses 
that of Proposition \ref{prop2c}).

Carefully examining the proof of Theorem \ref{thm3b}  
one verifies that there is no need to ever consider functions 
outside $\GG$ in case
this algebra is in the domain of the generators 
$\LL$, $\LL^*$, $\wt{\LL}$, $\LL^{\de f}_1$ and is closed under the
action of the corresponding Markovian semi-groups on $\CC_b$.  
It then suffices to define the operator $(\bA_1)_f$ of
(\ref{eq:aalpha}) only on $\GG$ and  
since here $\LL-\wt{\LL}=\bY$ we deduce from
(\ref{eq:gmdiff}) that  
$(\bA_1)_f g = \sum_{i=1}^d (\bX_i f) (\bX_i g) - f \bY g$
for all $g \in \GG$, as stated. Further, it follows that 
in such a situation $X^f_1(\cdot)$ has the response function 
corresponding to $(\bA_1)_f$ as soon as 
$\|\wt{\LL} g -\LL g\|^2 \leq - K \langle g \LL g \rangle_\mus$ 
for some $K<\infty$ and all $g \in \GG$, which is 
exactly our condition (\ref{eq:njd}).
\qed

Consider next the non-compact manifold $M=\Bbb R^d$, denoting 
by $\CC^\infty$, $\CC^\infty_b$, $\CC^\infty_0$  
the collections of smooth functions, smooth functions 
with bounded derivatives of all orders, smooth functions of 
compact support on $M=\Bbb R^d$, respectively. Let 
$\GG$ be the vector space spanned by $\won$ and the 
collection of Schwartz test functions on $M$ 
(i.e. functions in $\GG$ are 
elements of $\CC^\infty$ whose derivatives of all 
positive orders decay faster than any power of $\|x\|$),
and  
consider the Markovian semi-group $\bP_t$ on $\GG$ that is 
generated by $\LL=\sum_{i=1}^d \bX_i \circ \bX_i + \bX_0$ 
where for $i=0,\ldots,d$, 
$$
\bX_i =\sum_{k=1}^d a_{i,k}\frac{\partial}{\partial x_k} \,,
$$
with $a_{i,k}\in \CC^\infty_b$ for $i \ge 1$ while  
$a_{0,k}\in \CC^\infty$ 
with $\frac{\partial a_{0,k}}{\partial x_j} \in \CC^\infty_b$ 
for $k,j=1,\ldots,d$ (that is, the drift of our diffusion 
may be unbounded, but its derivatives are bounded). Hence, 
$$
\LL = \sum_{j,k=1}^d c_{j,k} \frac{\partial^2}{\partial x_j\partial x_k}  
+\sum_{j=1}^d b_j\frac{\partial}{\partial x_j} \,, 
$$
where for each $x \in M$, 
$$c_{j,k}=\sum_{i=1}^d a_{i,j} a_{i,k},
\quad b_j=a_{0,j}+\sum_{i,k=1}^d a_{i,k}\frac{\partial a_{i,j}}
{\partial x_k} \,.
$$
Suppose that $\{\bX_1,\ldots,\bX_d\}$ 
satisfies H\"ormander's strong hypo-elliptic 
condition and there exists a bounded below 
Lyapunov function $V\in \CC^\infty$, such that $\LL V\le 0$
and
$$\lim_{\|x\|\to\infty}V(x)=\infty.$$
This implies that the diffusion has 
a unique invariant measure $\mus$ with a
strictly positive smooth density $\rho\in \CC^\infty_b$
with respect to Lebesgue's measure on $M$
(see \cite{durrett}). 
In view of \cite[Theorem 3.14]{KuSt} 
the semi-group $\bP_t$ maps $\GG$ into itself (more precisely
they require bounded drift, but under our assumptions,  
for $\bX_0$ of linear growth the transition
probability function has sub-Gaussian tails
and once this is shown, a localization argument 
reduces to the case covered in \cite{KuSt}.
An alternative approach is to use that fact that if
$g\in\GG$ then $\LL \bP_t g =\bP_t (\LL g)$, where $\LL g\in\GG$ 
and use the weighted 
Sobolev norm estimates of \cite[Theorem 4.1]{EcHa}). 

Now take $f\in \CC^\infty_0$ and upon making the 
relevant modifications,
apply Proposition \ref{prop:diff} in this setting.

\section{Stochastic spin systems}\label{sec:spin}

We consider next systems of 
locally interacting diffusion processes, 
indexed by the $d$-dimensional lattice $\Bbb Z^d$.
Such processes naturally arise in statistical physics, 
where all Gibbs states of the interaction potential 
are invariant measures for the chosen dynamics.
In particular, in the presence of a 
phase transition we typically have non-uniqueness of 
the Gibbs state and infinitely many
invariant measures for the Markov process $X(\cdot)$. 
Note that in contrast with the setting of
Section \ref{sec:diff}, here the state space $\sss$
is such that we typically do not 
have an obvious dense algebra of test functions $\GG$ in
$\CCb$ which 
is closed under the action of the semi-group $\bP_t$
(for example, the algebra of functions depending on finitely
many coordinates is typically not closed under action of $\bP_t$).

For simplicity we restrict ourselves to 
pair interaction potentials and consider first the simpler 
case of compact spin spaces
$\sss= (\bS^1)^{\Bbb Z^d}$,
with spins taking values in
the one dimensional torus $\bS^1$ (equipped with 
Lebesgue measure and its $\sigma$-algebra $\sss_1$),
and having smooth, symmetric, finite range interactions. 
Specifically, for 
$\bx = (x_i, i \in \Bbb Z^d)$ and 
any $V\subset\subset\Bbb Z^d$ (i.e. $V$ finite),
consider the Hamiltonian
\begin{equation}\label{eq:poten}
H_V(\bx)=\sum_{i\in V} 
\Phi_{i}(x_i)+
\frac12\sum_{i \neq j \in V} \Phi_{i,j}(x_i,x_j)
+\sum_{i\in V, j\notin V} \Phi_{i,j}(x_i,x_j)\,,
\end{equation}
where  the potentials $\Phi_i\in \CC^\infty(\bS^1)$ and 
$\Phi_{i,j}\in \CC^\infty(\bS^1\times\bS^1)$, $i\ne j$ are 
such that $\Phi_{i,j}(x,y)=\Phi_{j,i}(x,y)=\Phi_{i,j}(y,x)$ 
and $\Phi_{i,j}
= {\bf 0}$ if $|i-j|>r$.
Let 
$$
H_i(\bx):=H_{\{i\}}(\bx) = \Phi_{i}(x_i)+
\sum_{j\in N(i)} \Phi_{i,j}(x_i,x_j)\,,
$$ 
where $N(i)=\{j\in \Bbb Z^d: 1<|i-j| \leq r\}$ denotes the 
$r$-neighborhood of $i$, excluding $i$.
Given smooth functions $\Psi_{i}(\bx)=\Psi_i(x_j, j \in N(i)) 
\in \CC^\infty((\bS^1)^{N(i)})$, set
$$b_i(\bx)=-\Phi_i'(x_i)-\sum_{j\in N(i)}\Phi'_{i,j}(x_i,x_j) 
+\Psi_{i}(\bx)e^{H_i(\bx)},$$
where 
$\Phi_i'(x)=\partial_x\Phi_i(x)$, 
$\Phi'_{i,j}(x,y)=\partial_x\Phi_{i,j}(x,y)$
(and hereafter $\partial_x$ denotes the smooth section 
corresponding to this differential operator on
the $\CC^\infty$-manifold $\bS^1$).
We assume that 
\begin{equation}
\sup_i|\Phi'_i(0)|\le C,\quad\sup_{i,j}|\Phi'_{i,j}(0,0)|\le C,
\label{eq:tag4}
\end{equation}
\begin{equation}
\inf_i\inf_{x}\Phi_i^{''}(x)\ge
 -K, \quad\inf_{i,j}\inf_{x,y}\Phi_{i,j}^{''}(x,y)\ge -K 
\label{eq:tag5}
\end{equation}
\begin{equation}
\label{eq:tag6}
\sup_{x,y} |\partial_x\partial_y\Phi_{i,j}(x,y)|\le c(i-j),
\end{equation}
\begin{equation}
\sup_i|\Psi_i({\bf 0})|\le C,\quad 
\sup_i\sup_{j\in N(i)}\sup_{\bx} |\partial_{x_j}\Psi_i(\bx)|\le C
\label{eq:tag7}
\end{equation}
for some finite constants $C$, $K$ and $c(k)$. 
Under these conditions there exists a 
strong Markov process $X(t)$ on $\sss=(\bS^1)^{\Bbb Z^d}$ 
(equipped with the corresponding product topology),
which is the unique strong solution of the
system of stochastic differential equations 
\begin{equation}
dX_i(t)= b_i(X(t))\,dt +\sqrt 2 dW_i(t),\qquad i\in\Bbb Z^d,
\label{eq:tag1}
\end{equation}
where $W_i(t)$, $i\in\Bbb Z^d$ are independent 
Brownian motions on $\bS^1$ (c.f. \cite{SS}).

We note in passing that with $\bS^1$ compact,
conditions (\ref{eq:tag4})-(\ref{eq:tag7}) 
trivially hold for example whenever $i \mapsto \Phi_i$, $i \mapsto \Psi_i$
and $(i,j) \mapsto \Phi_{i,j}$
are translation invariant (and also in some other cases). 

Let $\GG$ denote the dense subset of $\CC_b(\sss)$, 
consisting of all local smooth functions. That is, 
$\GG=\{\CC^\infty((\bS^1)^V )$ 
for some $V\subset\subset\Bbb Z^d\}$. Restricted to the 
algebra $\GG$,
the generator of the corresponding semi-group $\bP_t$ on $\CC_b(\sss)$
takes the form
\begin{equation}
\LL =\sum_i\big( \bX_i \circ \bX_i + b_i \bX_i \big) \,,
\label{eq:tag2}
\end{equation}
where $\bX_i g = \partial_{x_i}g$. Hence, 
$\LL g \in \GG$ for each $g \in \GG$. 

The solution $X(t)$ of (\ref{eq:tag1}) is smooth with respect
to $X(0)=\bx$ with derivatives 
$Y_{ij}(t)=\partial_{x_j} X_i(t)$ such that 
$Y_{ij}(t)=\sum_k \int_0^t B_{ik}(X(s)) Y_{kj}(s) ds$ and 
$Y_{ij}(0)=1_{j}(i)$, where by our assumptions 
$B_{ik} = \bX_k b_i$ are uniformly bounded and 
the sum is over the finite support of the
local function $b_i$. 
While $\bP_t g$ is not necessarily local for all $g \in \GG$ 
(so $\GG$ is not closed under the action of this semi-group),
nevertheless by the chain rule 
$\bX_i \bP_t g = \sum_k \E_{\bx} [ \bX_k g (X(t)) Y_{ik}(t) ]$,
which is thus uniformly bounded in $\bx$ (c.f. 
the proof of \cite[Theorem 2.2]{De-St2} or 
\cite{HS} which further shows the existence of 
a smooth transition probability density).
Further, iterating this procedure to deal with the second derivatives,
one easily verifies that 
for all $\wh{g}\in\wh{\GG}$,
\begin{equation}
\sum_i \|\bX_i \circ \bX_i \wh{g}\|_\infty+
\sum_i \|b_i \bX_i \wh{g} \|_\infty< \infty\,.
\label{eq:tag3}
\end{equation}
In particular, $h \bP_t g \in \DD(\LL)$ for all $h,g \in \GG$ and
the form (\ref{eq:tag2}) of $\LL$ extends to any such function.

Under our assumptions the collection 
${\bf G}(\Phi)$ of Gibbs measures 
associated with the potentials $\{\Phi_i, \Phi_{i,j}, i,j \in \Bbb Z^d \}$ 
is non-empty (see \cite{HS}). That is,
$\mu \in {\bf G}(\Phi)$ if and only if 
the DLR equations $\int_B \mu_V(d\by|\bx) \mu(d\bx) = \mu(B)$ 
hold for all $B \in \sss_1^{\Bbb Z^d}$ 
and $V \subset \subset \Bbb Z^d$, where  
the probability measure
$\mu_V(\cdot|\bx)$ on $(\bS^1)^V$ has the density 
$$
\mu_V(\by|\bx)=\frac{1}{Z_V(\bx)}\exp(-H_V(\by))
\prod_{i\notin V}\delta_{x_i}(y_i),
$$
with respect to Lebesgue measure on this set (and
$Z_V(\bx)\in\Bbb R_+$ is the corresponding 
normalizing constant). 
Since $\bS^1$ has no boundary and $\partial_{x_i} \Psi_i = 
{\bf 0}$,  
considering the form (\ref{eq:tag2}) and the
DLR relation for $V=\{i\}$, it is not hard to verify that 
$\langle \LL g \rangle_\mu =0$ for any
$\mu \in {\bf G}(\Phi)$ and $g \in \wh{\GG}$.
That is, any measure in ${\bf G}(\Phi)$ is invariant for 
$X(t)$. Further, 
fixing $\mus \in {\bf G}(\Phi)$, upon
using the explicit form of the local specification
$\mu_{\{i\}}(\,\cdot\,|\bx)$
and integration by part, it is easy to check that
for all 
$g,h \in \DD(\bX_i)$ and for $\mus$-a.e. $\bx \in \sss$,
\begin{equation}\label{eq:newam1}
\langle h \bX_i g \rangle_{\mu_{\{i\}}(\cdot|\bx)} =
\langle g (H'_{i} h  -\bX_i h) \rangle_{\mu_{\{i\}}(\cdot|\bx)}  
\end{equation}
where 
$H'_{i} = \bX_i H_i =\Phi_i'(x_i)+\sum_{j\in N(i)}\Phi'_{i,j}(x_i,x_j)$
is in $\GG$.
That is, 
the operator $\bX^{\star}_i=-\bX_i +H'_{i}$, having the same domain
as $\bX_i$, is the adjoint of $\bX_i$ in $L_2(\mus)$. 
In particular,  
$$
\oo{\LL} =\sum_i \big(-\bX^\star_i\circ \bX_i +\Psi_{i}e^{H_i}\bX_i\big) \,,
$$
with its adjoint in $L_2(\mus)$
$$
\oo{\LL}^\star  =\sum_i \big(-\bX^\star_i\circ \bX_i
 -\Psi_{i} e^{H_i}\bX_i\big)
$$
having the same domain as $\oo{\LL}$ and such that 
$\oo{\LL}^\star \GG \subseteq \GG$. Consequently, the same
applies for the self-adjoint operator  
$\wt{\LL} = - \sum_i \bX^\star_i\circ \bX_i$
(whose domain thus contains that of $\oo{\LL}$).

When $\Psi_i
={\bf 0}$ for all $i$, 
the Markov process $X(t)$ is merely 
the usual reversible Langevin dynamic 
and ${\bf G}(\Phi)$ is then precisely the collection of 
measures for which its semi-group $\bP_t$ is symmetric 
(see \cite[Theorem 4.3]{SS}). In this case,
it is further known that for $d=1,2$ there are no
other invariant measures for $\bP_t$ (c.f. \cite{HS}). 

Fixing a non-negative, smooth, local function 
$f\in \CC^\infty((\bS^1)^U)$ 
for some $U\subset\subset\Bbb Z^d$
(with $f 
\neq {\bf 0}$),
and $\mus\in{\bf G}(\Phi)$, the
probability measure $\bar\mu^f=Z_f^{-1} e^f \mus$ 
(with $Z_f=\langle e^f\rangle_\mus$) is a 
Gibbs measure corresponding to potentials
$\{\Phi_i, \Phi_{i,j}, i,j\in\Bbb Z^d$ and $\wt{\Phi}_U=-f\}$.
Consequently, $\bar \mu^f$ is invariant for 
the strong Markov process $X^f(t)$ 
which is the unique strong solution of
the SDE (\ref{eq:tag1}) with drift 
$$
b^f_i=-H_i'+\bX_i f + \Psi_i^f \exp(H_i-f 1_U(i)) \,,
$$
provided $\Psi^f_i \in \CC^\infty((\bS^1)^{N(i)})$
satisfy (\ref{eq:tag7}). For example, when taking 
$\Psi_i^f 
= {\bf 0}$ the generator of the corresponding
semi-group coincides with $\wt{\LL}^f + e^{-f} \wt{\LL}$
of Lemma \ref{lem2b} and Corollary \ref{cor2b}, 
at least when acting on smooth functions satisfying (\ref{eq:tag3}) 
(and in particular, on $\wh{\GG}$). 

We now seek $\Psi_i^f$ for which this Markovian perturbation 
has a response function. To this end, suppose first that 
$\Psi_i 
= {\bf 0}$ 
for all but finitely many $i \in \Bbb Z^d$. 
In this case
we may and shall enlarge $U$ so $\Psi_i 
= {\bf 0}$ 
for all $i \notin U$ 
in which case taking $\Psi_i^f = \Psi_i$ results 
effectively with the Markovian perturbation $X_1^f(\cdot)$ of 
Proposition \ref{prop2c}. 
The same argument we used to derive (\ref{eq:newam1}) also shows that 
$\langle r \Psi_i e^{H_i} \bX_i r \rangle_{\mu_{\{i\}}(\cdot|\bx)} = 0$
for $\mus$-a.e. $\bx \in \sss$ and any $r \in \DD(\bX_i)$.  
Consequently, for any $r \in \DD(\oo{\LL})$, 
$$
\langle r \oo{\LL} r \rangle_{\mus} = \langle r \wt{\LL} r \rangle_{\mus} = 
- \sum_i \|\bX_i r \|^2 \,,
$$
and applying Cauchy-Schwarz inequality we see that 
the condition (\ref{eq:sector}) 
holds here for $K = \sum_{i\in U}\|\Psi_i e^{H_i}\|^2$ finite. 
In view of Remark \ref{rem:inconv} we can apply Theorem \ref{thm3b} 
and conclude that this Markovian perturbation yields a response 
function (in the sense of Definition \ref{def:resp}), 
that corresponds to
\begin{equation}\label{eq:afdef}
\bA_f =\sum_i (\bX_i f) \bX_i - f \sum_{i \in U} \Psi_i e^{H_i} \bX_i \,.
\end{equation}
Whereas the condition (\ref{eq:sector}) typically fails to hold 
when $\Psi_i$ is non-zero at infinitely many sites $i \in \Bbb Z^d$, 
the perturbed drift $b_i^f$ with $\Psi_i^f = \Psi_i$ 
(where now $U=U_f$ remains the domain of $f$), can be directly 
shown to still yield the response function corresponding 
to the operator $\bA_f$ of (\ref{eq:afdef}). While we do not detail
this argument, note that we merely  
replaced the non-symmetric part of the perturbation 
(namely, $(1-e^{-f}) (\oo{\LL} r - \wt{\LL} r)$), by 
the ``localized'' non-symmetric operator
$(1-e^{-f}) \sum_{i \in U} \Psi_i e^{H_i} \bX_i r$,
which as we have seen, is dominated by 
$- \langle r \oo{\LL} r \rangle_{\mus}$ (and we can handle 
the latter as in the proof of Theorem \ref{thm3b}). 

Considering next the case of unbounded spins, where $\bS^1$
is replaced by $\Bbb R$, we restrict the state space 
$\sss$ to the subset of tempered configurations
$\bx\in \Bbb R^{\Bbb Z^d}$ such that $\sum_i|x_i|^2(1+|i|)^{-2p}$ is finite
for $p$ large enough. Considering only the usual Langevin dynamics
where $\Psi_i
= {\bf 0}$ 
and assuming once more
that (\ref{eq:tag4})--(\ref{eq:tag6}) hold, guarantees the
existence of a strong Markov process $X(t)$ with state space
$\sss$ which is the unique strong solution 
of the SDE (\ref{eq:tag1}), 
see \cite{SS}. 
In this setting we take 
$$\GG=\{g\in \CC^\infty(\Bbb R^V), V\subset\subset \Bbb Z^d, \partial_{x_i}g
\text{ has a compact support in } \Bbb R^V, \forall i\in V\}\,,$$
where again (\ref{eq:tag2}) holds and 
$\LL \GG \subseteq \GG$. 

Generally $\bP_t g$ is neither a local function, nor having  
derivatives of compact support, so once again 
$\GG$ is not closed under the action of $\bP_t$. 
Nevertheless, under suitable assumptions (\ref{eq:tag3}) holds,
the process $X(t)$ has invariant (Gibbs) measures and the
$\bar\mu^f$-symmetric Markovian semi-group 
associated with the perturbed drift $b_i^f$ yields the 
same explicit response function as before. For example,  
as shown in \cite{SS}, the set ${\bf G}(\Phi)$ of Gibbs measures is 
non-empty when for some $b>\sum_k c(k)$ 
$$
\sup_{x,i} ( b x^2 - x \Phi'_i(x) ) < \infty \,.
$$

\smallskip
\noindent {\bf Acknowledgment} 
We are grateful to Leticia Cugliandolo 
whose notes on the derivation of the FDT 
for reversible Markov chains on a
finite state space motivated this work and
to Yves Le Jan for his help with the proof of 
Proposition \ref{prop1}.
We also benefited from valuable feedback of
Jorge Kurchan,  
Joseph Avron, Joel Lebowitz and Andrea Montanari 
about the role and use of the FDT in physics, 
from a comment of Ben Goldys regarding the setting
for a non-compact state space, and from the anonymous 
referees input which greatly improved the
presentation of our work.

\end{document}